\documentclass[a4paper,11pt]{article}
\usepackage[usenames]{color}
\usepackage[utf8x]{inputenc}
\usepackage{amsmath}
\usepackage{amssymb}
\usepackage{amsfonts}
\usepackage{amsthm,amscd}
\newtheorem{theorem}{Theorem}[section]
\newtheorem{lemma}{Lemma}[section]
\newtheorem{corollary}{Corollary}[section]
\newtheorem{proposition}{Proposition}[section]
\theoremstyle{definition}
\newtheorem{definition}{Definition}[section]
\newtheorem{example}{Example}[section]

\theoremstyle{remark}
\newtheorem{remark}[theorem]{Remark}

\numberwithin{equation}{section}

\title{Hofer-like geometry and flux theory}
\author{\scshape 
S. Tchuiaga \thanks{tchuiagas@gmail.com, \ Department of Mathematics,  
The University of Buea, South West Region, Cameroon} }

\definecolor{couleurliens}{rgb}{1.0,0.,0.} 
 \definecolor{couleurliensref}{rgb}{0.,0.,1.} 
\definecolor{couleurliensurl}{rgb}{.3,.4,.3} 
\usepackage[pdftex,colorlinks=true,urlcolor=couleurliensurl,linkcolor=couleurliens,citecolor=couleurliensref]{hyperref}

\begin{document}
\maketitle\large 
\begin{abstract} 
This paper meticulously revisit and study the flux geometry of any 
compact connected oriented manifold $(M, \Omega)$. We 
generalize several well-known factorization results,
 exhibit some orbital conditions under which flux geometry can be studied, give a 
proof of 
the discreteness of the flux group for volume-preserving diffeomorphisms, deriving that any smooth isotopy in the kernel 
of the flux for volumr-preserving diffeomorphisms is a vanishing-flux path, and 
show that the kernel of the flux for volume-preserving diffeomorphisms is
 $C^1-$closed inside the group of all volume-preserving diffeomorphisms 
isotopic to the identity map: This recovers several well-known results from symplectic geometry.
The fix-points theory does not resist to the above machinery: We 
prove a general contractibility result with respect to the orbits of 
the fix-points for volume-preserving diffeomorphisms 
isotopic to the identity map via vanishing-flux paths, generalize and    
solve the Arnold conjecture using the Thurston fragmentation property.  
In the sequel, we use fix-points to: Characterize the flux 
geometry of certain  
$C^0-$limits of sequences of vanishing-flux paths and  volume-preserving diffeomorphisms. Beside 
this, a $C^0-$criterion for the existence of at least one fix-point is given, and 
a weak version of the generalized $C^0-$flux conjecture is solved. Finally, 
we construct a pseudo right-invariant metric on the group of all volume-preserving diffeomorphisms 
isotopic to the identity map,   
 prove several comparison results suitable to the study of the Hofer-like geometry of the group $Ham(N,\omega)$, 
of all Hamiltonian diffeomorphisms of a closed symplectic manifold $(N,\omega)$,
 derive the equivalence between  
the Hofer and the Hofer-like metrics on $Ham(N,\omega)$, and exhibit a computational proof 
of the non-degeneracy of the Hofer-like energies: Here, an
 outcome is that the Calabi group  
controls the Hofer-like 
geometry of the group $Ham(N,\omega)$ of any non-simply connected closed symplectic manifold $(N,\omega)$. 
This includes several other interesting results. \\ 
\end{abstract}

\section{Introduction}\label{SC00}
On a compact manifold $M$, given any closed differential $p-$form 
 $\alpha$  together with an isotopy $\Phi = (\phi_t) $, they satisfy the following equation   
\begin{equation}\label{0F1}
 (\phi_t)^\ast\alpha - \alpha = d\left( \mathcal{F}_\alpha(\Phi)(t)\right)
\end{equation}
for all $t$, where $
 \mathcal{F}_\alpha
(\Phi)(t) := \int_0^{t}(\Phi_s)^\ast(\alpha(\dot\Phi_s))ds,$ 
for all $t$, and $d$ is the usual differential operator such that $d\circ d = 0$. 
Equation (\ref{0F1}) is involved in the study of the de Rham cohomology of 
the manifold $M$: 
It had been used by Moser \cite{MS}, Banyaga \cite{Ban78, Ban97}, Calabi \cite{Cal} and others in the 
study of the geometry of the group of volume-preserving 
diffeomorphisms (e.g. it follows from  Banyaga \cite{Ban78}, 
Thurston \cite{TH}, and Weinstein \cite{Wien} works that 
if two given volume-preserving 
isotopies are homotopic relatively to fix endpoints, then they  
have the same flux). So far, this is the most known (not to say the only) condition 
which guaranties the equality of 
the fluxes of two volume-preserving 
isotopies. What seems to be a bit  ambiguous here is that, the flux of 
each isotopy is a de Rham cohomology class in some $H^\ast(M,\mathbb{R})$ ( 
the $\ast-$th de Rham group of $M$ with real coefficients), while the group $H^\ast(M,\mathbb{R})$ 
only depend on the topology of $M$, and does not depend on any differentiable structure 
on $M$.  So, one might expect that analyzing the orbits (in $M$)  of a volume-preserving 
isotopy of $(M, \Omega)$ (where $\Omega$ is the orientation form on $M$) may ''determine'' the 
flux of the latter isotopy. Of course, 
the difficulty to face (or the prize to pay) here is the large number of the orbits 
generated by a single isotopy on a manifold: It looks almost impossible to describe 
all of them. But, this does not tell us that 
such a study is not possible. The optimal cost here could reside in the way to look at the situation.
 The question addressed  in this paper is that of the existence 
of other(s) criterion (or criteria)  under which one can study the flux geometry of 
a compact oriented manifold.\\ 
We shall see that investigating the above criteria will not only enrich the theory of 
flux geometry, but also the dynamics of fix-points for volume-preserving diffeomorphisms, and 
the theory of Hofer-like geometry of 
a closed 
oriented manifold. To do so, we organize the paper as follows:\\ 
In Section \ref{Prelim} we introduce the definitions 
of some useful tools that we shall need. The third section
is devoted to revisits the construction of the flux homomorphism of a 
compact connected oriented manifold, and generalize 
several factorization results found in \cite{Ban97,TD2}: 
The concerned results are,  
Lemma \ref{Al-c2}, Lemma \ref{Al-c3}, 
Proposition \ref{SC41}, Proposition \ref{Fact3}, and Proposition \ref{Important}.\\   
 In Section \ref{Displa}, we study the displacements of closed $1-$forms 
(Proposition \ref{Impor}, and Proposition \ref{GF1}). This includes: 
  Theorem \ref{OrTheo0} which  
gives orbital conditions under which the fluxes of two isotopies with the same endpoints are equal, 
Theorem \ref{OrTheo3} which states that the flux group for volume-preserving diffeomorphisms
 is a discrete group,  Lemma \ref{Rigiditylem0} which states that the sub-group of 
vanishing-flux volume-preserving diffeomorphisms 
is $C^1-$closed inside the group of all volume-preserving diffeomorphisms isotopic to the identity map, 
and 
Lemma \ref{Rigiditylem2}  which states that the sub-group of all vanishing-flux 
volume-preserving isotopies is $C^0-$closed in the group of all 
volume-preserving isotopies.
The fifth section deals with fix-points theory. Here we prove:
Theorem \ref{Functor3} shows 
that the fix-points of vanishing flux volume-preserving diffeomorphisms are exactly the 
vanishing points of a certain smooth function real valued on $M$ modulo the space of closed $1-$forms, 
Lemma \ref{FixPt4} which states that the orbit of 
any fix-point of the time-one map of any vanishing-flux volume-preserving isotopy 
is contractible, Theorem \ref{Functor11} which generalizes and solves the 
Arnold conjecture, and Theorem \ref{Fix-Points3-5-7} gives a $C^0-$condition that 
predicts the minimum number of fix-points of 
a vanishing-flux volume-preserving diffeomorphisms. This includes, Lemma \ref{Rigiditylem2} and 
Lemma \ref{Rigiditylem3} which show that fix-points can determine the flux geometry of 
the $C^0-$limit of a sequence of vanishing-flux volume-preserving diffeomorphisms,
 Theorem \ref{Rigidity} which states  that if $\gamma$ is  
the $C^0-$limit of a sequence of vanishing flux volume-preserving isotopies (loops), then 
 any of smooth orbit of $\gamma$ is contractible. Also, Lemma \ref{Separation} separates 
the orbits of a non-vanishing 
flux volume-preserving isotopy from minimal 
geodesics between their endpoints. Section \ref{Hofer-like} is devoted to the study of the metric 
geometry of the group of volume-preserving 
diffeomorphisms isotopic to the identity map. Here, 
we construct a right-invariant pseudo-metric on the latter group (Proposition \ref{Fixpoints-metric2}), 
and study the Hofer-like geometry with 
the intension to highlight how the latter is related to flux geometry.   
We prove several comparison results illustrating the impacts of flux geometry in the  study of the Hofer-like 
geometry of the group $Ham(M,\omega)$: The concerned results are,
 Lemma \ref{Vanishingamma}, Lemma \ref{Nonvanishingamma}, 
Proposition \ref{Positivelength}, Proposition \ref{proEQU}, and Theorem \ref{EQU}). 
In the sequel, a proof of the non-degeneracy of a suitable class of Hofer-like energies follows, 
and a computational proof 
of the equivalence between the Hofer norm and 
the restriction of the Hofer-like norm to the group $Ham(M,\omega)$ is given: A remarkable 
 outcome here is that the Calabi group controls the Hofer-like geometry of $Ham(M,\omega)$. 
\section{Preliminaries}\label{Prelim}
Let $M$ be an $n-$dimensional closed Riemannian manifold equipped with an orientation form $\Omega$.
 A diffeomorphism $\phi:M\rightarrow M$ preserves a differential $p-$form  
$\alpha$ on $M$ if $\phi^\ast(\alpha) = \alpha$. 
Given any differential $p-$form  $\alpha$ on $M$, let
 $Diff^\infty(M,\alpha)$ denote the group of all diffeomorphisms from $M$ to $M $ that preserve  $\alpha$. 
\subsection{Isotopies}
An isotopy $\Phi = (\phi_t)$ of $M$ is a smooth map from $[0,1]$ into  $Diff^\infty(M)$ such that $\phi_0 = id_M$ 
 (see \cite{Ban78, BanTch1, McDuff-SAl} for more details).  Let $Diff_0^\infty(M)$ denote the set of 
all time-one maps of smooth isotopies in $Diff^\infty(M)$, let $Iso(M, \alpha)$ 
denote the space of all smooth isotopies in
 $Diff^\infty(M,\alpha)$, and by $G_\alpha(M)$, we denote the set of all  
time-one maps of all isotopies of $Iso(M, \alpha)$.  
For each $ \psi \in G_\alpha(M)$, we shall denote by 
  $Iso(\psi)_\alpha$  the set $\{\Phi\in Iso(M,\alpha)| \Phi(1) = \psi\}$. 
Any  smooth isotopy $\Phi = (\phi_t)$ on $M$ gives rise 
to a smooth family of smooth vector fields  $(\dot\phi_t)$ over $M$ defined by  
$
 \dot\phi_t(\phi_t(x)) = \dfrac{d}{dt}(\phi_t(x)),$
for all $t$, and for all $x\in M$. 
\subsubsection{Volume-preserving isotopies}
An isotopy $\Phi = (\phi_t)$ on $(M, \Omega)$ is said to be volume-preserving if 
its associated smooth family of smooth vector fields  $(\dot\phi_t)$ consists of 
free divergence vector fields i.e. for each $t$
 the $(n-1)-$form $ \imath(\dot\phi_t)\Omega$ is closed. This is equivalent 
to require that $\Phi\in Iso(M, \Omega)$. 
\subsubsection{Concatenation operation for isotopies}\label{Concatenation}
It is always possible to construct a smooth increasing function $f:[0,1]\rightarrow[0,1] $ such that its 
restriction to an interval $[0,\delta] $ vanishes while its restriction to $[(1-\delta),1]$ is the 
constant function $1$ 
with $ 0 \textless\delta\leq \dfrac{1}{8} $.
Assume this done, and consider the smooth functions $ \lambda(t) = f(2t)$ for all $ 0\leq t\leq \frac{1}{2}$, and 
$\tau(t) = f(2t -1)$ for all  $ \frac{1}{2}\leq t\leq 1.$ If
 $\Phi = (\phi_t)$ and  $\Psi = (\psi_t)$ are two given isotopies, then one defines an 
isotopy $(\Phi\ast_r\Psi)$ with time-one 
map $\phi_1\circ\psi_1$ by setting: 
$$
(\Phi\ast_r\Psi)(t) =
\left\{
\begin{array}{l}
\phi_{\lambda(t)},  \hspace{0,2cm} if \hspace{0,2cm} 0\leq t\leq \frac{1}{2},\\
\phi_1\circ\psi_{\tau(t)},  \hspace{0,2cm} if \hspace{0,2cm}  \frac{1}{2}\leq t\leq 1.    \\
\end{array}
\right.
$$
The isotopy $(\Phi\ast_r\Psi)$ constructed above is called the 
right concatenation of the isotopy $\Phi$ by the isotopy $\Psi$. Similarly, we 
can define the left concatenation of the isotopy $\Phi$ by the isotopy $\Psi$, denoted $(\Psi\ast_l\Phi)$, as 
map $\psi_1\circ\phi_1$ by setting: 
$$
(\Psi\ast_l\Phi)(t) =
\left\{
\begin{array}{l}
\phi_{\lambda(t)},  \hspace{0,2cm} if \hspace{0,2cm} 0\leq t\leq \frac{1}{2},\\
\psi_{\tau(t)}\circ\phi_1,  \hspace{0,2cm} if \hspace{0,2cm}  \frac{1}{2}\leq t\leq 1.   \\
\end{array}
\right.
$$
An orbit of a point $p\in M$ under $(\Phi\ast_r\Psi)$ 
(resp. under $(\Psi\ast_l\Phi)$) is the orbit of $p$ under 
$\Phi$ ''glued'' with the image under $\phi_1$ of the orbit of $p$ under $\Psi$ (resp. 
is the orbit of $p$ under $\Phi$ ''glued'' with the orbit of $\phi_1(p)$ under $\Psi$). 
Through all the paper, for each isotopy $\Phi$ and each point $p\in M$, we shall denote by
 $\mathcal{O}_p^\Phi$ the orbit of $x$ under the action of $\Phi$, 
while $-\mathcal{O}_p^\Phi$ will represent the orbit 
of $p$ under the action of $\Phi^{-1}.$
\subsection{The de Rham groups}\label{De Rham}
Let $H^\ast(M,\mathbb{R})$ denote the $\ast-$th de Rham cohomology group (with real coefficients) of $M$, 
and let $\mathcal{Z}^\ast(M)$ denote the space of closed $\ast-$forms on $M$. Fix a linear section of the natural 
projection $P: \mathcal{Z}^\ast(M)\rightarrow H^\ast(M,\mathbb{R}),$ and denote 
it as $\mathcal{S}: H^\ast(M,\mathbb{R})\rightarrow \mathcal{Z}^\ast(M)$. Each 
$\alpha \in \mathcal{Z}^\ast(M)$ decomposes as $$
 \alpha = \mathcal{S}(P(\alpha)) + (\alpha - \mathcal{S}(P(\alpha)),$$
with $P(\mathcal{S}(P(\alpha))) = P(\alpha)$ and $P((\alpha - \mathcal{S}(P(\alpha))) = 0$. 
Note that a form $\alpha$ such that $P(\alpha) = 0$ is called an exact form. 
So, in the decomposition 
above, for each $\alpha \in \mathcal{Z}^\ast(M)$, the form $\mathcal{S}(P(\alpha))$ 
is called the non-exact part of $\alpha$ 
while $ (\alpha - \mathcal{S}(P(\alpha))$ is 
called the exact part of $\alpha$. 
Now, let $\mathcal{H}^\ast(M, \mathcal{S})$ denote the 
image $\mathcal{S}(P(\mathcal{Z}^\ast(M))\subset \mathcal{Z}^\ast(M)$: 
$\mathcal{H}^\ast(M, \mathcal{S})$ is isomorphic to $H^\ast(M,\mathbb{R})$, and 
it follows from Hodge's theory that the latter set is  
a finite dimensional vector space over $\mathbb{R}$ whose dimension is denoted 
$b_\ast(M),$ is called the 
$\ast-$th Betti number of the manifold $M$ \cite{FWa}. 
For $\ast = 1$, we 
fix a basis ${(H_i)_{1\leq i\leq b_1(M)}}$ on $\mathcal{H}^1(M, \mathcal{S})$ and equip it with a norm $|.|$
 (see \cite{TD2} for more details).\\ 
We shall also write $\overline{\mathcal{B}^\ast_{ \mathcal{S}}(0,1)} $ to mean the closure
of the open unit ball $\mathcal{B}^\ast_{ \mathcal{S}}(0,1)$
of the finite dimensional vector space $\mathcal{H}^\ast(M, \mathcal{S})$: The set 
$\overline{\mathcal{B}^\ast _{\mathcal{S}}(0,1)}$ 
is compact. 
\section{The flux homomorphism revisited}\label{Flux} 
 Here, we 
 revisit the construction of the flux homomorphism and derive some consequences. We will 
need the following formula: For each closed $1-$form $\alpha$, we have 
\begin{equation}\label{Al-o1}
 \mathcal{F}_{\alpha}(\Phi\circ\Psi)(t) = \mathcal{F}_{\alpha}(\Psi)(t) 
+ \mathcal{F}_{\alpha}(\Phi)(t)\circ\Psi(t) + \delta_t(\Phi,\Psi,\alpha),
\end{equation}
where $\Phi$  and $ \Psi$ are two isotopies of a closed manifold, and $\delta_t(\Phi,\Psi,\alpha)$ is a 
constant which depends on $\Phi,$ $ \Psi,$ $\alpha$ and $t$. 
\begin{lemma}(\cite{TD2})\label{Al-c1}$ $ Let $(M, \omega)$ be a closed symplectic manifold, and let $\alpha$ be a closed $1-$form. If $\Phi$ and $ \Psi$ are two symplectic isotopies,\\ then we have  
$
 \delta_t(\Phi,\Psi,\alpha) = 0,
$
for all $t$. 
\end{lemma}

According to  Lemme \ref{Al-c1}, on a closed symplectic manifold, relation (\ref{Al-o1}) becomes $
 \mathcal{F}_{\alpha}(\Phi\circ\Psi)(t) = \mathcal{F}_{\alpha}(\Psi)(t) + \mathcal{F}_{\alpha}(\Phi)(t)\circ\Psi(t), 
$
for all $t$. The latter equality generalize to volume-preserving case, and then seems to tell us that for each fixed closed $1-$form $\alpha$, 
the map $\Phi\mapsto \int_M\mathcal{F}_{\alpha}(\Phi)(1)\Omega,$ 
induces a group homomorphism from the space $Iso(M,\Omega)$ into $\mathbb{R}$.\\

 We have the following facts. 
\begin{lemma}\label{Al-c2} $ $ For each fixed closed $1-$form $\alpha$, and for each volume-preserving isotopy 
 $\Phi$, the  integral  $\int_M\mathcal{F}_{\alpha}(\Phi)(1)\Omega$ 
is independent of the choice of any representative $\beta$
in the de Rham cohomology class  $ [\alpha]$. 
\end{lemma}
\begin{proof} $ $
  Let $\beta\in [\alpha],$ i.e.  
 $\alpha - \beta = df$ for some smooth function $f:M\rightarrow\mathbb{R}$. Assume that  
$\Phi = \{\phi_t\},$ and derive from a direct computation that,
 $$0 = \iota(\dot\phi_t)[(\alpha -\beta)\wedge\Omega]
= \iota(\dot\phi_t)[(\alpha -\beta)]\Omega 
+ (\alpha -\beta)\wedge\iota(\dot\phi_t)\Omega.$$ 
Since  the $(n-1)-$form $\imath(\dot\phi_t)\Omega$ is closed and $\alpha -\beta = df$, then 
$ (\alpha -\beta)\wedge\iota(\dot\phi_t)\Omega = d(f\iota(\dot\phi_t)\Omega).$ Hence, we have 
$ - (\iota(\dot\phi_t)\alpha)\Omega  +(\iota(\dot\phi_t)\beta)\Omega =  d(f\iota(\dot\phi_t)\Omega),$
and integrating the above relation over $M$ yields, 
$ - \int_M(\iota(\dot\phi_t)\alpha)\Omega  + \int_M(\iota(\dot\phi_t)\beta)\Omega 
=  \int_Md(f\iota(\dot\phi_t)\Omega),$
where  the right hand side vanishes because of Stokes' Theorem since $M$ is closed.
\end{proof}

\begin{lemma}\label{Al-c3} $ $
If $\Phi$ and $ \Psi$ are two isotopies which are 
homotopic relatively to fix endpoints, then 
$ \mathcal{F}_{\alpha}(\Phi)(1) = \mathcal{F}_{\alpha}(\Psi)(1),$
for each closed  $1-$form $\alpha$. 
\end{lemma}
\begin{proof}$ $
  Let $\alpha$ be a closed $1-$form; set  
$ \Psi = (\psi_t),$ and  
$\Phi = (\phi_t).$ We have to prove that
$\mathcal{F}_{\alpha}(\Phi)(1)(x) = \mathcal{F}_{\alpha}(\Psi)(1)(x),$ 
for all $x\in M$. To that end, note that  
a homotopy between $ \Psi$ and 
$\Phi$ is a smooth map $H: I\times I\rightarrow Diff^\infty_0(M)$ such that
$H(0,t) = \psi_t$, $H(1,t) =\phi_t$ for all $t$; $H(s, 0) = id_M$ and $H(s,1)= \phi_1 =\psi_1$ for all $s$. 
Thus, for each $x\in M$, the homotopy $H$ induces 
an homotopy $\overline{H}_x$ between the orbits $\mathcal{O}^\Phi_x:t\mapsto\phi_t(x)$ 
and $\mathcal{O}^\Psi_x:t\mapsto\psi_t(x)$; i.e.   
$\mathcal{O}^\Phi_x$ and $\mathcal{O}^\Psi_x$ constitute the boundary of a $2-$chain in $M$. 
Therefore, Stokes' Theorem implies that 
$\int_{\mathcal{O}^\Phi_x}\alpha =\int_{\mathcal{O}^\Psi_x}\alpha,$ 
because $\alpha$ is closed. To conclude, observe that for each $x\in M$, we also have
$\mathcal{F}_{\alpha}(\Phi)(1)(x) = \int_{\mathcal{O}^\Phi_x}\alpha,~$
and $\mathcal{F}_{\alpha}(\Psi)(1)(x) = \int_{\mathcal{O}^\Psi_x}\alpha$. 
\end{proof}

Lemme \ref{Al-c1} and Lemme \ref{Al-c2} suggest that: 
Each volume-preserving isotopy   
$\Phi$ induces an element $\widetilde{S}_\Omega(\Phi)\in 
Hom(H^{1}(M,\mathbb{R}),\mathbb{R})\cong H^{n-1}(M,\mathbb{R}),$ defined as follows: 
To any $[\alpha]\in H^{1}(M,\mathbb{R}) $ one assigns  
$\widetilde{S}_\Omega(\Phi)([\alpha]) := \int_M\mathcal{F}_{\alpha}(\Phi)(1)\Omega.$
Thus, we have a group homomorphism 
$\widetilde{S}_\Omega: Iso(M,\Omega)\rightarrow H^{n-1}(M,\mathbb{R}),$
which is surjective. This motivated the following factorization result that generalizes Proposition $2.4-$\cite{TD2}.  
\begin{proposition}\label{SC41}{\bf (Factorization I)} $ $ 
  Let $(M,\Omega)$ be a closed oriented manifold. 
Let $\alpha = (\alpha_t)$ be a 
smooth family of smooth closed $1-$forms, and let 
$ \Phi = \{\phi_t\}\in Iso(M, \Omega)$.  Then, for each $t\in [0,1],$ we have, 
$
 \int_M\mathcal{F}_{\alpha_t}(\Phi)(t)\Omega = \langle  P(\alpha_t), \widetilde{S}_\Omega(\bar{\Phi}_t)\rangle,
$ 
 where $\bar{\Phi}_t$ is the isotopy $s\mapsto\phi_{st},$  and  
$\langle,\rangle: H^1(M,\mathbb{R})\times H^{n-1}(M,\mathbb{R})\rightarrow \mathbb{R}, 
([\alpha], [\beta])\mapsto \int_{M}\alpha\wedge\beta,$
is the usual Poincar\'e pairing. 
\end{proposition}
\begin{proof}$ $
 
Fix  an arbitrary $t$ in $[0,1]$, and since the form  
 $\alpha_t\wedge\Omega$ is of degree $(n + 1)$ on an $n-$dimensional manifold, then it is trivial. 
This implies that
$ \alpha_t(\dot\phi_s)\Omega  - \alpha_t \wedge\iota(\dot\phi_s)\Omega= 0,$
for each $s\in[0,t]$. Composing each member of the above equality by 
$\phi_s^\ast$ gives 
\begin{equation}\label{UB1}
 \phi_s^\ast\left(\alpha_t(\dot\phi_s)\right)\Omega - 
 \phi_s^\ast(\alpha_t)\wedge \phi_s^\ast\left(\iota(\dot\phi_s)\Omega\right) = 0,
\end{equation}
for each $s\in[0,t]$. Using (\ref{0F1}), one deduces that for each $s\in[0,t]$, we have 
\begin{equation}\label{UB01}
\phi_s^\ast(\alpha_t) = \alpha_t + df_{\{\phi_t\}, \alpha_t}^s,
\end{equation}
 where 
 $f_{\{\phi_t\}, \alpha_t}^s := \int_0^s\alpha_t(\dot\phi_u)\circ\phi_udu.$ 
 Observe that  (\ref{UB1}) and (\ref{UB01})  together implies that,
 \begin{equation}\label{UB2}
   \phi_s^\ast\left(\alpha_t(\dot\phi_s)\right)\Omega 
=  \alpha_t\wedge\phi_s^\ast\left(\iota(\dot\phi_s)\Omega\right)+ 
 df_{\{\phi_t\}, \alpha_t}^s\wedge\phi_s^\ast\left(\iota(\dot\phi_s)\Omega\right),
 \end{equation}
for each  $s\in[0,t]$. 
Thus, 
\begin{equation}\label{UB22}
 \int_{M}\left(\int_0^t\phi_s^\ast\left(\alpha_t(\dot\phi_s)\right)ds\right)\Omega = 
\int_{M}\alpha_t\wedge\left(\int_0^t \phi_s^\ast\left(\iota(\dot\phi_s)\Omega\right)ds\right)
\end{equation}
$$ + \int_{M}\left(\int_0^t\left(\phi_s^\ast\left(\iota(\dot\phi_s)\Omega\right)\wedge 
df_{\{\phi_t\}, \alpha_t}^s\right) ds\right).$$ 
On the other hand, since the form $\phi_s^\ast\left(\iota(\dot\phi_s)\Omega\right) $ is closed, then 
\begin{equation}\label{UB222}
 \int_{M}\int_0^t\phi_s^\ast\left(\iota(\dot\phi_s)\Omega\right)\wedge df_{\{\phi_t\}, \alpha_t}^sds 
=  \int_{M}d[\int_0^t f_{\{\phi_t\}, \alpha_t}^s\phi_s^\ast\left(\iota(\dot\phi_s)\Omega\right)ds],
\end{equation}
and since $M$ is without boundary, then Stokes' theorem implies  
\begin{equation}\label{UB2222}
 \int_{M}d[\int_0^t f_{\{\phi_t\}, \alpha_t}^s\phi_s^\ast\left(\iota(\dot\phi_s)\Omega\right)ds] 
= \int_{\partial M}\int_0^t f_{\{\phi_t\}, \alpha_t}^s\phi_s^\ast\left(\iota(\dot\phi_s)\Omega\right)ds = 0,
\end{equation}
for each $t$. Combining (\ref{UB22}), (\ref{UB222}), and (\ref{UB2222}) together yields
\begin{equation}\label{UB3}
 \int_{M}\left(\int_0^t\phi_s^\ast\left(\alpha_t(\dot\phi_s)\right)ds\right)\Omega = 
  \int_{M}\alpha_t\wedge\left(\int_0^t \phi_s^\ast\left(\iota(\dot\phi_s)\Omega\right)ds\right).
\end{equation}
For each $t$, the de Rham cohomology class of the  
$(n-1)-$form $\int_0^t \phi_s^\ast\left(\iota(\dot\phi_s)\Omega\right)ds$ is exactly 
 $\widetilde{S}_\Omega(\bar{\Phi}_t),$ i.e.  (\ref{UB3}) becomes 
$
 \int_{M} \mathcal{F}_{\alpha_t}(\Phi)(t)\Omega  = 
\langle P(\alpha_t),\widetilde{S}_\Omega(\bar{\Phi}_t) \rangle,
$
for all $t$. 
\end{proof}

Let us mention that for each closed $1-$form $\alpha$ (fixed), 
we have a surjective mapping $\Phi \mapsto\int_M\mathcal{F}_{\alpha}(\Phi)(1)\Omega$, 
 from $Iso(M,\Omega)$ onto $ H^{n-1}(M,\mathbb{R})$. Let $r$ be a real number: If $r= 0$, 
then take $\Phi$ to be the constant 
path identity. Suppose $r\neq 0$: This means that there exists 
at least $\Phi = \{\phi^t\}\in Iso(M,\Omega)$ such that 
$\int_M\mathcal{F}_{\alpha}(\Phi)(1)\Omega\neq 0$. With Proposition \ref{SC41}, this is equivalent 
to say that $$ \int_0^1\langle P(\alpha), P(\iota(\dot\phi^t)\Omega) \rangle dt \neq 0,$$ i.e., 
there exists a $t_0\in [0, 1]$ such that $ \langle P(\alpha), P(\iota(X)\Omega)\rangle \neq 0$, with  
$X :=  \dot\phi^{t_0}$. 
Let $\Psi$ be the flow generated by $X$, and define another conservative 
vector field by setting: $Y:= rX/\langle P(\alpha), P(\iota(X)\Omega)\rangle$. If $\Theta$ stands 
for the flow generated by $Y$, then  one has:
$$\int_M\mathcal{F}_{\alpha}(\Theta)(1)\Omega = \langle P(\alpha),\widetilde{S}_\Omega(\Theta) \rangle =
 \langle P(\alpha),P(\iota(Y)\Omega)\rangle $$ 
$$= r\langle P(\alpha),P(\iota(X)\Omega)\rangle /\langle P(\alpha), P(\iota(X)\Omega)\rangle = r.$$  
Furthermore, for each non-exact closed $1-$form $\alpha$ (fixed), since the 
linear map $L_\alpha: H^{(n-1)}(M,\mathbb{R})\ni\beta \mapsto \langle P(\alpha),\beta \rangle$, is non-trivial, 
hence the latter is surjective. Now, if we denote by $H_\alpha$ the surjective mapping 
$ \Phi\mapsto \int_M\mathcal{F}_{\alpha}(\Phi)(1)\Omega$, then by Proposition \ref{SC41} 
we have $H_\alpha = L_\alpha \circ \widetilde{S}_\Omega$, which is surjective, hence so is 
the mapping  
$\widetilde{S}_\Omega$. On the 
other hand, let $\{\Phi_i = (\phi^t_i)\}\subset Iso(M,\Omega)$ be a sequence that converges 
in the $C^0-$topology to $\Phi= \{\phi^t\}\in Iso(M,\Omega)$, i.e., 
$\Phi_i\circ \Phi^{-1}\xrightarrow{C^0}Id$. Then, it is not hard derive with the help of Lemma 3.10-\cite{TD2}
 that, 
$$\lim_i|\int_M\mathcal{F}_{\alpha}(\Phi_i\circ \Phi^{-1})(1)\Omega|
\leq 4|\alpha|_0Vol_{\Omega}(M)d_{C^0}(\Phi_i\circ \Phi^{-1}, Id)\Omega\rightarrow0, i\rightarrow\infty.
 $$
Thus, by Proposition \ref{SC41}, this means that for each closed $1-$form $\alpha$, 
we have\\ $ \lim_i\langle P(\alpha), \widetilde{S}_\Omega(\Phi_i\circ \Phi^{-1})\rangle = 0$, 
i.e. 
$ \langle P(\alpha), \lim_i\widetilde{S}_\Omega(\Phi_i\circ \Phi^{-1})\rangle = 0$, for all 
 closed $1-$form $\alpha$. Hence, $\lim_i\widetilde{S}_\Omega(\Phi_i)= \widetilde{S}_\Omega(\Phi)$: The 
map $\widetilde{S}_\Omega$ is a continuous group 
homomorphism. $\maltese$\\

The  
homomorphism $\widetilde{S}_\Omega$ had been studied by Weinstein \cite{Wien} and 
Thurston \cite{TH} when $\Omega$ is a volume form. 
Set $\varGamma_\Omega=\widetilde{S}_\Omega(\pi_1(G_\Omega(M)))$, and 
 pass to the quotient 
to obtain the following commutative diagram: 
\begin{equation}\label{diagram}
 \begin{array}{ccc}
\stackrel{\sim}{G_\Omega}\!(M) & \stackrel{\widetilde{S}_\Omega}{\longrightarrow} & H^{(n-1)}(M,\mathbb{R}) \\ 
P_1\downarrow &     & \downarrow P_2 
\\ G_\Omega(M) &  \stackrel{S_\Omega}{\longrightarrow} & H^{(n-1)}(M,\mathbb{R})/\varGamma_\Omega, 
\end{array}
\end{equation}
where  $\stackrel{\sim}{G_\Omega}\!(M)$ represents the quotient space of $Iso(M, \Omega)$ 
with respect to the equivalent 
relation ''homotopic relatively to fixed endpoints''; $P_i$ are projection maps. 
The following fact are well-known:
\begin{itemize}
 \item The sub-group $\varGamma_\Omega$ 
is discrete, (unpublished result of Thurston): We shall give an proof of this result later on using arguments from global analysis.
\item The group $\ker S_\Omega$ has the fragmentation property, \cite{TH}.
\item As a consequence of Moser' theorem, the group $ G_\Omega(M)$  is locally 
connected by smooth arcs, \cite{Ban97}.
\item Boothby proved that the group $ G_\Omega(M)$  is $p-$transitive, \cite{Boo}.
\item An interesting question is that: Is any smooth path in $\ker S_\Omega$ a vanishing-flux path? We shall answr this question later on. 
\end{itemize}
Proposition \ref{SC41}  gives a short proof of a generalized Lemme \ref{Al-c1}: Infact, 
(\ref{Al-o1}) implies that\\ 
$$\int_{M}\mathcal{F}_{\alpha}(\Phi\circ\Psi)(t)\Omega =
 \int_{M}\mathcal{F}_{\alpha}(\Psi)(t)\Omega + \int_{M}\mathcal{F}_{\alpha}(\Phi)(t)\circ\Psi(t)\Omega
 + \delta_t(\Phi,\Psi,\alpha)\int_{M}\Omega,$$ 
for all $t$, and according to Proposition \ref{SC41}, 
we have $$\int_{M}\mathcal{F}_{\alpha}(\Phi\circ\Psi)(t)\Omega 
= \int_{M}\mathcal{F}_{\alpha}(\Psi)(t)\Omega + \int_{M}\mathcal{F}_{\alpha}(\Phi)(t)\Omega.$$
So, combining the above equalities, it follows that $\delta_t(\Phi,\Psi,\alpha)\int_{M}\Omega =0, $ for all $t$; 
this implies that $\delta_t(\Phi,\Psi,\alpha) = 0,$ for all $t$. $\maltese$\\

Here are some consequences of Proposition \ref{SC41}. 
\begin{lemma}\label{Orlem1}$ $  
 Let $\Phi\in Iso(M,\Omega)$ be a loop at the identity map. For each closed $1-$form
$\alpha,$ the smooth function $x\mapsto\int_{\mathcal{O}_{x}^\Phi}\alpha$ is constant 
and agrees with  
$\dfrac{1}{Vol_\Omega(M)}\langle P(\alpha),\widetilde{S}_\Omega(\Phi)\rangle$. 
\end{lemma}
\begin{proof} $ $
 Let $\alpha$ be a  closed $1-$form. If $P(\alpha) = 0$, then 
$\int_{\mathcal{O}_{x}^\Phi}\alpha = 0,$ for all $x\in M$ since $\mathcal{O}_{x}^\Phi$ is a closed curve. 
On the other hand,  
suppose $P(\alpha) \neq 0$, and set $\mathcal{H} = \mathcal{S}(P(\alpha))$. 
Since $\Phi(0) = \Phi(1) = id_M$, by formula (\ref{0F1}) the map 
$x\mapsto\int_{\mathcal{O}_{x}^\Phi}\alpha = \mathcal{F}_{\alpha}(\Psi)(1)(x)$ is constant. 
Thus, Proposition \ref{SC41} implies that
$$
 Vol_\Omega(M)\int_{\mathcal{O}_{x}^\Phi}\alpha = Vol_\Omega(M)\int_{\mathcal{O}_{x}^\Phi}\mathcal{H} = 
\int_M\mathcal{F}_{\alpha}(\Psi)(1) \Omega = 
\langle P(\mathcal{H}),\widetilde{S}_\Omega(\Phi)\rangle,
$$
for all $x \in M$.
\end{proof}
\begin{remark}\label{Obs1}
 Geometrically, Lemma \ref{Orlem1} seems to suggest that on a closed oriented manifold  
$(M,\Omega)$, the set 
$\pi_1(G_\Omega(M))\backslash (\ker\widetilde{S}_\Omega\cap\pi_1(G_\Omega(M))$ measures the 
obstruction that impeaches the orbits generated by all the elements of $\pi_1(G_\Omega(M))$ 
to be all contractible in $M$. $\maltese$ 
\end{remark}

\begin{lemma}\label{Orlem2} $ $ Let $(M,\Omega)$ be a closed oriented manifold. 
 Let $\Phi\in Iso(M,\Omega)$ be a 
loop at the identity.  
Then, $\widetilde{S}_\Omega(\Phi) = 0$,  if and only if, there exists 
 $z\in M$ such that the orbit $\mathcal{O}_{z}^\Phi$ is contractible. 
\end{lemma}
\begin{proof} $ $
  As in the proof of Lemma \ref{Orlem1}, we have 
$$\dfrac{1}{Vol_\Omega(M)}\langle P(\alpha), \widetilde{S}_\Omega(\Phi)\rangle 
= \int_{\mathcal{O}_{x}^\Phi}\alpha,$$
for each closed $1-$form $\alpha$: If  $z\in M$ is 
such that the orbit $\mathcal{O}_{z}^\Phi$ is contractible, 
then Stokes' theorem implies that  $ \int_{\mathcal{O}_{z}^\Phi}\alpha = 0$  
for each closed $1-$form $\alpha$; and this 
imposes that 
$ \langle P(\alpha), \widetilde{S}_\Omega(\Phi)\rangle  = 0,$ for all 
closed $1-$form $\alpha$, i.e. 
$\widetilde{S}_\Omega(\Phi) = 0$. 
Conversely, if $\widetilde{S}_\Omega(\Phi)=0,$ 
then we immediately see that $ 0= \int_{\mathcal{O}_{z}^\Phi}\alpha,$ 
 for all 
closed $1-$form $\alpha$, and for all $z\in M$. This implies that the orbit
 $\mathcal{O}_{z}^\Phi $ must be contractible.
\end{proof}

\begin{proposition}\label{Fact3}{\bf (Factorization II)} $ $
Let $(M,\Omega)$ be a compact connected oriented manifold with boundary.
  Let $\alpha = (\alpha_t)$ be a 
smooth family of smooth closed $1-$forms, and let 
$ \Phi = \{\phi_t\}\in Iso(M, \Omega)$.  Then, for each $t\in [0,1],$ we have,\\
$$
 \int_M\mathcal{F}_{\alpha_t}(\Phi)(t)\Omega = \langle  P(\alpha_t), \widetilde{S}_\Omega(\bar{\Phi}_t)\rangle 
+ \int_{\partial M}\left( \int_0^t f_{\Phi, \alpha_t}^s\phi_s^\ast\left(\iota(\dot\phi_s)\Omega\right)ds\right),
$$ 
where $\bar{\Phi}_t$ is the isotopy $s\mapsto\phi_{st}$, and  
 $f_{\Phi, \alpha_t}^s := \int_0^s\alpha_t(\dot\phi_u)\circ\phi_udu$. 
\end{proposition}
\begin{remark}\label{Obs2} $ $ From Proposition \ref{Fact3}, we can argue
 that on a compact connected oriented manifold 
$(M,\Omega)$ with boundary $\partial M$, for any loop
 $\Phi =\{\phi_t\}\in (\ker\widetilde{S}_\Omega\cap\pi_1(G_\Omega(M))$, we have  
$$
 \int_{\mathcal{O}_{z}^\Phi}\alpha = \frac{1}{Vol_\Omega(M)}\int_{\partial M}
\left(\int_0^1  f_{\Phi, \alpha}^s\phi_s^\ast\left(\iota(\dot\phi_s)\Omega\right)ds\right),
$$
 for each $z\in M$, and each closed $1-$form $\alpha$ on $M$, 
where 
 $f_{\Phi, \alpha}^s := \int_0^s\alpha(\dot\phi_u)\circ\phi_udu$: Geometrically, 
this seems to suggest that 
on a compact connected oriented manifold $(M,\Omega)$
with boundary, the smooth manifold $\partial M$ could be 
 an obstruction that impeaches certain orbits for volume-preserving isotopies on $M$ to be contractible. $\maltese$ 
\end{remark}

\begin{proposition}\label{Important}{\bf (Factorization III)} $ $ Let $(M,\Omega)$ 
be a closed oriented manifold 
possibly with boundary 
such that $\Omega = \alpha^l := \underbrace{\alpha\wedge\dots\wedge\alpha}_{l-times},$ where $\alpha$ 
is a closed $\deg(\alpha)-$form with
$\deg(\alpha)= \dim(M)/ l$, and for each volume-preserving isotopy 
$\Phi = \{\phi_t\}
$ the $[\deg(\alpha)-1]-$form  $ \iota(\dot\phi_t)\alpha$ 
is closed for each $t$. Therefore, 
\begin{enumerate}
 \item if  $\deg(\alpha)$ is even, then\\  $\widetilde{S}_\Omega(\Phi) = \dfrac{\dim(M)}{\deg(\alpha)}
P  \left(
\alpha^{\left(\left(\dfrac{\dim(M)}{\deg(\alpha)} \right)-1\right)}\wedge
\int_0^1\iota(\dot\phi_t)\alpha dt \right),$ 
\item if  $\deg(\alpha)$ is odd, then\\  $\widetilde{S}_\Omega(\Phi) 
= (-1)^{\left(1 + E\left( \dfrac{\dim(M)}{2\deg(\alpha)}\right)\right)}
P\left(
\alpha^{\left(\left(\dfrac{\dim(M)}{\deg(\alpha)} \right)-1\right)}\wedge
\int_0^1\iota(\dot\phi_t)\alpha dt\right),$
\end{enumerate}
 where for each positive real number $x$, $E(x)$ is x minus its fractional part; 
and  $P$ is the natural 
projection $P: \mathcal{Z}^\ast(M)\rightarrow H^\ast(M,\mathbb{R}).$
\end{proposition}
\begin{proof} $ $
 Assume $\alpha$ is a closed $\deg(\alpha)-$form of even degree, and let $\beta$ be any 
closed $1-$form. Since $\beta\wedge\Omega = 0,$ 
we derive as in the proof of Proposition \ref{SC41} that 
\begin{equation}\label{F1UB1}
 \phi_s^\ast\left(\beta(\dot\phi_s)\right)\Omega =
\phi_s^\ast(\beta\wedge \iota(\dot\phi_s) \Omega) =  \beta\wedge\phi_s^\ast( \iota(\dot\phi_s) \Omega) 
+ df^s_{\Phi, \beta}\wedge\phi_s^\ast( \iota(\dot\phi_s) \Omega) 
\end{equation}
$$ = l\beta\wedge \phi_s^\ast\left(\iota(\dot\phi_s)\alpha\right)\wedge\phi_s^\ast\left(\alpha^{l-1}\right) 
+ df^s_{\Phi, \beta}\wedge \phi_s^\ast( \iota(\dot\phi_s) \Omega) $$
for each $s\in[0,1]$, with $ l =\dim(M)/\deg(\alpha)$, i.e. 
\begin{equation}\label{F2UB3}
 \int_{M}\left(\int_0^1\phi_s^\ast\left(\beta(\dot\phi_s)\right)ds\right)\Omega = 
  l\int_{M}\beta\wedge\left(\int_0^1 \phi_s^\ast\left(\iota(\dot\phi_s)\alpha\wedge\alpha^{l-1}\right)ds\right) 
\end{equation}
$$ + \int_{\partial M}\int_0^1
\left( f_{\Phi, \alpha_t}^s\phi_s^\ast\left(\iota(\dot\phi_s)\Omega\right)ds\right).$$
Thus, we derive by the mean of Proposition \ref{Fact3} that,
$$\langle  P(\beta), \widetilde{S}_\Omega(\Phi)\rangle  = 
\int_{M}\left(\int_0^1\phi_s^\ast\left(\beta(\dot\phi_s)\right)ds\right)\Omega = \langle P(\beta), 
lP\left(\alpha^{l-1}\wedge\int_0^1 \iota(\dot\phi_s)\alpha ds\right)\rangle,$$
for all closed $1-$form $\beta$. Therefore, the non-degeneracy of $\langle,\rangle$ implies that\\ 
$ \widetilde{S}_\Omega(\Phi) = 
lP\left(\alpha^{l-1}\wedge\int_0^1 \iota(\dot\phi_s)\alpha ds\right).$
 For the second item, assume that $\alpha$ is a closed $\deg(\alpha)-$form of odd degree, 
and let $\beta$ be any closed $1-$form. Similar calculations lead to 
\begin{equation}\label{F3UB1}
 \phi_s^\ast\left(\beta(\dot\phi_s)\right)\Omega = (-1)^\epsilon 
 \phi_s^\ast(\beta)\wedge \phi_s^\ast\left(\iota(\dot\phi_s)\alpha\right)\wedge\phi_s^\ast\left(\alpha^{l-1}\right) 
+ df^s_{\Phi, \beta}\wedge\phi_s^\ast( \iota(\dot\phi_s) \Omega),
\end{equation}
for each $s\in[0,1]$, with 
$\epsilon = \left(1 + E\left(\dfrac{\dim(M)}{2\deg(\alpha)} \right)\right),$  
where for each positive real number $x$, $E(x)$ is x minus its fractional part. 
\end{proof}

\begin{remark}
If $(M,\omega)$ is a $2l-$dimensional symplectic manifold oriented with the Liouville volume 
form $\Omega_0 = \dfrac{\omega^l}{l!}$, then  
set $\alpha = \omega$, and derive from 
 Proposition \ref{Important} that 
$$\widetilde{S}_{\Omega_0}(\Phi) = \dfrac{l}{l!}
P(\int_0^1\phi^\ast_t(\iota(\dot\phi_t)\omega)dt\wedge\omega^{(l-1)}) 
= \dfrac{1}{(l-1)!}\widetilde{S}_{\omega}(\Phi)\wedge P(\omega^{(l-1)}) ,$$
for all symplectic isotopy $\Phi$:  
Proposition \ref{Important} is a generalization of the factorization result found by Banyaga \cite{Ban78, Ban97}.
\end{remark}
\section{The geometry of the displacements of closed $1-$forms}\label{Displa}
In this section, we assume that $\partial M =\emptyset$. 
Let $\psi\in Diff^\infty(M)$,
 for each fixed $p\in M$ (such a point is called a base point), and  
for all $z\in M$, let $\xi_z$ denote a geodesic from  $p$ to $z$. Then, for each closed $1-$form  
$\alpha,$ the following smooth function $
 \nu^{\psi, \alpha}_{p} : M \rightarrow \mathbb{R}, z\mapsto \int_{\xi_z}\left(\psi^\ast\alpha - \alpha\right),
$
is well-defined since it does not depend on the choice of any curve from $p$ to $z$. 
The real number $\int_{\xi}\left(\psi^\ast\alpha - \alpha\right)$  
can be viewed in a certain sense as the displacement of $\alpha$ under 
the action of $\psi$ along the curve $\xi$.\\

We have the following fact.

\begin{corollary}\label{Impor}{\bf(Changing base points)} $ $ Let $\psi\in Diff^\infty_0(M)$
and $\alpha$ be any closed non-trivial $1-$form. If $\xi$ and $\gamma$  
are two curves such that $\xi(1) = \gamma(1)$, 
then we have\\ $\nu^{\psi, \alpha}_{\xi(0)} - 
\nu^{\psi, \alpha}_{\gamma(0)} + \nu^{\psi, \alpha}_{\gamma(0)}(\xi(0)) = 0.$  
\end{corollary}
\begin{proof} $ $ 
 
 Since the curves $\xi$, $\gamma$, and $C$ form the boundary of a $2-$chain in $M$, and by (\ref{0F1}) 
the $1-$form  $\psi^\ast\alpha - \alpha$ 
is exact, then 
$
 \int_{\xi}\left(\psi^\ast\alpha - \alpha\right)  - \int_{\gamma}\left(\psi^\ast\alpha - \alpha\right)
 + \int_{C}\left(\psi^\ast\alpha - \alpha\right) = 0,
$
i.e.  $\nu^{\psi, \alpha}_{\xi(0)} - \nu^{\psi, \alpha}_{\gamma(0)} + \nu^{\psi, \alpha}_{\gamma(0)}(\xi(0)) = 0.$  

\end{proof}

For further investigations, let us consider the following  function defined on $ Diff^\infty_0(M)$  
 as follows: For each $p\in M$ (fixed), then to each $ \psi \in Diff^\infty_0(M)$, 
assign the quantity $ \varDelta(\psi, \alpha)_p$ 
defined by:
\begin{equation}\label{Functor2}
 \varDelta(\psi, \alpha)_p :=
\left\{
\begin{array}{l}
 0, \hspace{0.1cm} \hspace{0.1cm}if\hspace{0.1cm} \alpha = 0\\
 \dfrac{1}{\|\alpha\|_{L^2}}\int_M\nu^{\psi, \alpha}_{p}\Omega 
\hspace{0.1cm} \hspace{0.1cm},  \hspace{0.1cm}if\hspace{0.1cm} 
\alpha\in\left(\mathcal{Z}^1(M)\backslash\{0\}\right),
\end{array}
\right.  
\end{equation}
where $\|\alpha\|_{L^2}$ is the $L^2-$Hodge norm of $\alpha$:   $\|\alpha\|_{L^2}^2:=\int_M \alpha\wedge\ast\alpha$, 
where $\ast$ is the usual Hodge star operator.   
The above formula is well-defined for each $p\in M$ and all $ \alpha\in \mathcal{Z}^1(M)$. 

\begin{corollary}\label{GF10}   $ $ Fix a point $x\in M$, and an element 
$\alpha\in\left(\mathcal{Z}^1(M)\backslash\{0\}\right)$. 
For each $\psi \in G_\Omega(M),$ we have, 
$$
\varDelta(\psi, \alpha)_x = 
\dfrac{1}{\|\alpha\|_{L^2}}
\langle P(\alpha),\widetilde{S}_\Omega(\Psi) \rangle-
\dfrac{Vol_\Omega(M)}{\|\alpha\|_{L^2}}\int_{\mathcal{O}_{x}^\Psi}\alpha,
$$
and the real number 
 $\varDelta(\psi, \alpha)_x$ does not depend on any choice of an isotopy in  $Iso(\psi)_\Omega$.
\end{corollary}
\begin{proof} $ $
 Fix a point  $x\in M$.  
Formula (\ref{0F1}) implies that if
$\Phi,\Psi\in Iso(\psi)_\Omega$, then 
$ \mathcal{F}_\alpha(\Psi)(1)(y) = \mathcal{F}_\alpha(\Phi)(1)(y) + cte,$
for all $y\in M$, i.e.
$$ \mathcal{F}_\alpha(\Psi)(1)(y) - \mathcal{F}_\alpha(\Psi)(1)(x) 
= \mathcal{F}_\alpha(\Phi)(1)(y) - \mathcal{F}_\alpha(\Phi)(1)(x),$$
for all $y\in M$.
On the other hand,  if we let $\xi_y$ denote any curve from  $x$ to $y$, 
then it follows that 
$\nu^{\psi, \alpha}_{x}(y)= \mathcal{F}_\alpha(\Psi)(1)(y) 
- \mathcal{F}_\alpha(\Psi)(1)(x),$ 
for all $y\in M$, i.e. 
$$ \mathcal{F}_\alpha(\Psi)(1)(y) - \mathcal{F}_\alpha(\Psi)(1)(x) = \nu^{\psi, \alpha}_{x}(y) 
= \mathcal{F}_\alpha(\Phi)(1)(y) - \mathcal{F}_\alpha(\Phi)(1)(x),$$
for all $y\in M$.
 Hence, integrating the above relation over $M$ yields, 
 $$ \int_M\mathcal{F}_\alpha(\Psi)(1)\Omega - \mathcal{F}_\alpha(\Psi)(1)(x)\int_M \Omega
 = \int_M\nu^{\psi, \alpha}_{x}\Omega 
= \int_M\mathcal{F}_\alpha(\Phi)(1)\Omega - \mathcal{F}_\alpha(\Phi)(1)(x)\int_M\Omega.$$
Therefore, multiplying the above equality in both sides by $ 1/\|\alpha\|_{L^2}$ yields
  $$\dfrac{1}{\|\alpha\|_{L^2}}\int_M\mathcal{F}_\alpha(\Psi)(1)\Omega   
- \dfrac{Vol_\Omega(M)}{\|\alpha\|_{L^2}}\mathcal{F}_\alpha(\Psi)(1)(x)
 = \dfrac{1}{\|\alpha\|_{L^2}}\int_M\nu^{\psi, \alpha}_{x}\Omega,$$ 
$$ 
=  \dfrac{1}{\|\alpha\|_{L^2}}\int_M\mathcal{F}_\alpha(\Phi)(1)\Omega  
- \dfrac{Vol_\Omega(M)}{\|\alpha\|_{L^2}}\mathcal{F}_\alpha(\Phi)(1)(x),$$
for all $\Phi,\Psi\in Iso(\psi)_\Omega$, where  $\mathcal{F}_\alpha(\Psi)(1)(x) 
= \int_{\mathcal{O}_{x}^\Psi}\alpha,$ 
while  Proposition \ref{SC41} implies that, 
$\int_{M} \mathcal{F}_{\alpha}(\Phi)(1)\Omega  = 
\langle P(\alpha),\widetilde{S}_\Omega(\Phi) \rangle.$
\end{proof}

Let $\Phi$ be an isotopy. For each positive integer $l$, we shall write $\Phi^l$  
to mean the $l-$fold of the isotopy $\Phi$, i.e. $\Phi^l = \underbrace{\Phi\circ\dots\circ\Phi}_{l-times}.$ 
Similarly, for each diffeomorphism $\phi$ and for each positive integer $l$, we shall write $\phi^l$  
to mean $\phi^l = \underbrace{\phi\circ\dots\circ\phi}_{l-times}$.  

\begin{proposition}\label{GF1}{\bf(Fragmentation)}$ $
 For each $x\in M$ (fixed), and for each $\alpha \in \mathcal{Z}^1(M)\backslash\{0\}$, 
the mapping
$
\varDelta(., \alpha)_x  : G_\Omega(M) \rightarrow \mathbb{R}, 
\psi \mapsto \varDelta(\psi, \alpha)_x,
$
has the following properties:
\begin{enumerate}

\item For all $\psi^1,\dots, \psi^k \in G_\Omega(M),$  we have 
$$
 \varDelta(\psi^1\circ\dots\circ\psi^k, \alpha)_x = \varDelta(\psi^k, \alpha)_x + 
\sum_{k-1\geq i\geq 1}\varDelta(\psi^i, \alpha)_{(\psi^{i+ 1}\circ\dots\circ\psi^k)(x)}.$$ 
\item In particular, 
$$\varDelta(\psi^1\circ\psi^2, \alpha)_x= \varDelta(\psi^1, \alpha)_x + \varDelta(\psi^2, \alpha)_x 
+ \dfrac{Vol_\Omega(M)}{\|\alpha\|_{L^2}}\left(\int_{\mathcal{O}_{x}^\Phi}\alpha- 
\int_{\mathcal{O}_{\psi^2(x)}^\Phi}\alpha\right),$$
 for all $\Phi\in Iso(\psi^1)_\Omega.$
\item  For all $l\in \mathbb{Z}^\ast$, and for each $\psi \in G_\Omega(M)$, we have 
$$
 \varDelta(\psi^l, \alpha)_x = l\varDelta(\psi, \alpha)_x + 
\dfrac{Vol_\Omega(M)}{\|\alpha\|_{L^2}}\left( l\int_{\mathcal{O}_x^{\Phi^\epsilon}}\alpha - 
\sum_{i = 0}^{\epsilon l -1}\int_{\mathcal{O}^{\Phi^\epsilon}_{\psi^{\epsilon i}(x)}}\alpha \right),
$$
\end{enumerate}
 where $\Phi\in Iso(\psi)_\Omega$, and $\epsilon = sign(l)$, and 
 $\psi^{\epsilon i} = \underbrace{\psi^\epsilon\circ\dots\circ\psi^\epsilon}_{i-times}.$
\end{proposition}
\begin{proof} $ $

For the first item, we proceed by induction: For $k =2$, we derive from Corollary \ref{GF10} that 
$$\varDelta(\psi^1\circ\psi^2, \alpha)_x =  \dfrac{1}{\|\alpha\|_{L^2}}\langle  
P(\alpha), \widetilde{S}_\Omega(\Psi\ast_l \Phi)\rangle - 
\dfrac{Vol_\Omega(M)}{\|\alpha\|_{L^2}}\int_{\mathcal{O}_{x}^{\Psi\ast_l \Phi}}\alpha,$$
where $\Phi\ast_l\Psi\in Iso(\psi^1\circ\psi^2)_\Omega$ is the left concatenation of $\Psi\in Iso(\psi^2)_\Omega$ by  
$\Phi\in Iso(\psi^1)_\Omega.$ Observe that 
$\widetilde{S}_\Omega(\Psi\ast_l\Phi) = \widetilde{S}_\Omega(\Psi) + \widetilde{S}_\Omega(\Phi);$ and 
$\mathcal{O}_{x}^{\Phi\ast_l\Psi} = \mathcal{O}_{x}^{\Psi}\bigsqcup \mathcal{O}_{\psi^2(x)}^{\Phi},$ 
and deduce that 
$\varDelta(\psi^1\circ\psi^2, \alpha)_x = \varDelta(\psi^2, \alpha)_x
+ \varDelta(\psi^1, \alpha)_{ \psi^2(x)}$. Let $k \textgreater 2$: Assume that 
$$\varDelta(\psi^1\circ\dots\circ\psi^{k-1}, \alpha)_x = \varDelta(\psi^{k-1}, \alpha)_x + 
\sum_{k-2\geq i\geq 1}\varDelta(\psi^i, \alpha)_{(\psi^{i + 1}\circ\dots\circ\psi^{k-1})(x)},$$
set $\phi:= \psi^1\circ\dots\circ\psi^{k-1}$, and compute 
 $$\varDelta(\phi\circ\psi^{k}, \alpha)_x = \varDelta(\psi^k, \alpha)_x +  
\varDelta(\phi, \alpha)_{\psi^k(x)} $$ 
$$= \varDelta(\psi^k, \alpha)_x + \varDelta(\psi^{k-1}, \alpha)_{(\psi^{k}(x))} + 
\sum_{k-2\geq i\geq 1}\varDelta(\psi^i, \alpha)_{(\psi^{i + 1}\circ\dots\circ\psi^{k-1})(x)}
, 
$$ 
$$= \varDelta(\psi^k, \alpha)_x + 
\sum_{k-1\geq i\geq 1}\varDelta(\psi^i, \alpha)_{(\psi^{i + 1}\circ\dots\circ\psi^{k-1})(x)}
.$$ 
For the second item, a direct computation 
with the help of Corollary \ref{GF10} 
shows that 
$$\varDelta(\psi^1, \alpha)_{ \psi^2(x)} =  \dfrac{1}{\|\alpha\|_{L^2}}\langle  
P(\alpha), \widetilde{S}_\Omega(\Phi)\rangle - 
\dfrac{Vol_\Omega(M)}{\|\alpha\|_{L^2}}\int_{\mathcal{O}_{\psi^2(x)}^\Phi}\alpha$$ 
$$= \varDelta(\psi^1, \alpha)_x + 
\dfrac{Vol_\Omega(M)}{\|\alpha\|_{L^2}}\left(\int_{\mathcal{O}_{x}^\Phi}\alpha- 
\int_{\mathcal{O}_{\psi^2(x)}^\Phi}\alpha\right).$$
For the third item, assume $l$ positive, and for each $i = 0,1,\dots ,l-1;$ 
consider a smooth increasing function  $\lambda_i$, 
defined from
$[i/l, i+1/l]$ onto $[0,1]$. Define an isotopy $\Phi^l$ as follows: $\Phi^l(t) = \Phi(\lambda_i(t))\circ\psi^i,$
whenever $t\in [i/l, i+1/l]$, and $i = 0,1\dots ,l-1$. Since $\Phi\in Iso(\psi)_\Omega$, 
we have $\Phi^l\in Iso(\psi^l)_\Omega$. Therefore, 
 we obtain 
 $\widetilde{S}_\Omega(\Phi^l) = l\widetilde{S}_\Omega(\Phi), \mbox{and} \hspace{0.2cm}
  \mathcal{O}_{x}^{\Phi^l} = \bigcup_{i = 0}^{l-1}\mathcal{O}_{\psi^i(x)}^{\Phi},$
and this in turn with Corollary \ref{GF10} implies that 
$$\varDelta(\psi^l, \alpha)_x = \dfrac{1}{\|\alpha\|_{L^2}}\langle  
P(\alpha), \widetilde{S}_\Omega(\Phi^l)\rangle - 
\dfrac{Vol_\Omega(M)}{\|\alpha\|_{L^2}}\int_{\mathcal{O}_{x}^{\Phi^l}}\alpha
 $$ $$= \dfrac{l}{\|\alpha\|_{L^2}}\langle  P(\alpha), \widetilde{S}_\Omega(\Phi)\rangle  - 
\dfrac{Vol_\Omega(M)}{\|\alpha\|_{L^2}}\sum_{i = 0}^{l-1}\int_{\mathcal{O}_{\psi^i(x)}^{\Phi}}\alpha$$
$$
 =l\varDelta(\psi, \alpha)_x +  \dfrac{Vol_\Omega(M)}{\|\alpha\|_{L^2}}
\left(l\int_{\mathcal{O}_{x}^{\Phi}}\alpha - 
\sum_{i = 0}^{l-1}\int_{\mathcal{O}_{\psi^i(x)}^{\Phi}}\alpha\right).$$
When $\epsilon = -1,$  i.e. $l$ is negative, the result follows from similar arguments.
\end{proof}

\begin{lemma}\label{Cont1}$ $  For each $x\in M$ (fixed),  
and for each $\alpha \in \mathcal{Z}^1(M)$, 
 the map\\ 
$\varDelta( . ,\alpha)_x :  Diff^\infty_0(M) \rightarrow \mathbb{R}, 
\psi \mapsto \varDelta(\psi, \alpha)_x,$
is continuous with respect to the $C^0-$topology on $ Diff^\infty_0(M) $.
\end{lemma}
\begin{proof} $ $
  Let $(\psi_i)\subseteq  Diff^\infty_0(M)$ be a sequence  which converges to 
$\psi \in  Diff^\infty_0(M)$ with respect to the $C^0-$metric. Suppose that 
for  $i$ sufficiently large $d_{C^0}(\psi_i,\psi) \leq r(g)/i$ where $r(g)$ is the injectivity radius of the 
Riemannian metric 
$g$ on $M$ (fixed). Since for all $i$ sufficiently large we have the metric condition, 
$d_{C^0}(\psi_i,\psi) \leq r(g)/i$ holds, 
then for all $y\in M$, 
the points $\psi_i(y)$ and $\psi(y)$  
can be connected through a unique minimal geodesic $\chi_y^i$.  Let $z$ be any point of $M$ which realizes the 
supremum of the function  $y\mapsto |\nu^{\psi_i,\alpha}_{x}(y) - \nu^{\psi,\alpha}_{x}(y)|,$ and for all  
$i$ sufficiently large, denote by    
 $\square(z,x,\psi_i, \psi)$ the smooth $2-$chain delimited by the curves $\psi_i(\gamma_z)$, 
$\psi(\gamma_z)$, $\chi_z^i$, and $\chi_{x}^i$, where $\gamma_z$ is any smooth curve $x$ to $z$. 
Compute, 
$$ \nu^{\psi_i,\alpha}_{x}(z) - \nu^{\psi,\alpha}_{x}(z) = \int_{\gamma_z}\psi_i^\ast(\alpha) 
- \int_{\gamma_z}\psi^\ast(\alpha)
 = \int_{\psi_i(\gamma_z)}\alpha - \int_{\psi(\gamma_z)}\alpha,$$
 and since  $\alpha$ is closed, derive from  Stokes' theorem that 
$0 = \int_{\square(z,x,\psi_i, \psi)} d\alpha = \int_{\partial\square(z,x,\psi_i, \psi)} \alpha,$
i.e. $$ \sup_y|\nu^{\psi_i,\alpha}_{x}(y) - \nu^{\psi,\alpha}_{x}(y)| 
= |\int_{\gamma_z}\psi_i^\ast(\alpha) - \int_{\gamma_z}\psi^\ast(\alpha)|
 = |\int_{\psi_i(\gamma_z)}\alpha - \int_{\psi(\gamma_z)}\alpha|$$
$$= |\int_{\chi_z^i}\alpha - \int_{\chi_{x}^i}\alpha|
 \leq|\int_{\chi_z^i}\alpha| +  |\int_{\chi_{x}^i}\alpha|
\leq 2\|\alpha\|_{L^2}d_{C^0}(\psi_i,\psi),$$
because the speed of a minimal geodesic is bounded from above by the distance between its endpoints. 
Hence, 
$$|\varDelta(\psi_i, \alpha)_x-\varDelta(\psi, \alpha)_x|\leq 
\dfrac{Vol_\Omega(M)}{\|\alpha\|_{L^2}}\sup_y|\nu^{\psi_i,\alpha}_{x}(y) - 
\nu^{\psi,\alpha}_{x}(y)|$$ $$\leq 2Vol_\Omega(M)d_{C^0}(\psi_i,\psi)\rightarrow0,i\rightarrow\infty.$$
\end{proof}
\subsection{Flux geometry via orbits}
\begin{theorem}\label{OrTheo0} $ $ Let $(M,\Omega)$ be a closed oriented manifold. 
For each $\varphi\in \left( G_\Omega(M)\backslash\{id_M\}\right),$  let 
$\Phi$ and $\Psi$ be two volume-preserving isotopies in $Iso(\varphi)_\Omega$.  
 If there exists a point  $z_0\in M$ 
for which the orbits $\mathcal{O}_{z_0}^{\Psi}$ and $\mathcal{O}_{z_0}^{\Phi}$ are 
homotopic relatively to fix endpoints, then\\ $\widetilde{S}_\Omega(\Phi) = \widetilde{S}_\Omega(\Psi)$. 
\end{theorem}
\begin{proof} $ $
 
Suppose that  the $1-$cycle 
$(-\mathcal{O}_{z_0}^{\Psi})\bigsqcup\mathcal{O}_{z_0}^{\Phi}$ is contractible. 
Consider $\alpha$ to be any closed $1-$form such that $P(\alpha) \neq 0$,
 and let $\mathcal{H} = \mathcal{S}(P(\alpha)).$
 Corollary \ref{GF10} implies that 
\begin{equation}\label{IFG}
 \varDelta(\varphi, \alpha)_{z_0} = \dfrac{1}{\|\alpha\|_{L^2}}\langle P(\alpha), 
\widetilde{S}_\Omega(\Theta)\rangle - 
\dfrac{Vol_\Omega(M)}{\|\alpha\|_{L^2}}\int_{\mathcal{O}_{z_0}^\Theta}\alpha,
\end{equation}
for all $\Theta\in Iso(\varphi)_\Omega$. In particular, since each 
of the isotopies $\Phi$ and $\Psi$ belongs to $Iso(\varphi)_\Omega,$  
we then obtain from (\ref{IFG}) that\\
$
 \langle P(\alpha), \widetilde{S}_\Omega(\Phi)\rangle  - 
Vol_\Omega(M)\int_{\mathcal{O}_{z_0}^\Phi}\alpha = 
 \langle p(\alpha), \widetilde{S}_\Omega(\Psi)\rangle- 
Vol_\Omega(M)\int_{\mathcal{O}_{z_0}^\Psi}\alpha,$\\
i.e., $
 \langle P(\alpha),\widetilde{S}_\Omega(\Psi)- \widetilde{S}_\Omega(\Phi)\rangle = 
 Vol_\Omega(M)\int_{(-\mathcal{O}_{z_0}^{\Psi})\bigsqcup\mathcal{O}_{z_0}^{\Phi}}\alpha,$
for all closed  $1-$form $\alpha$. But, 
since $ (-\mathcal{O}_{z_0}^{\Psi})\bigsqcup\mathcal{O}_{z_0}^{\Phi}$ is contractible,
and $\alpha$ is closed, one derives from Stokes' theorem that: $
 \int_{(-\mathcal{O}_{z_0}^{\Psi})\bigsqcup\mathcal{O}_{z_0}^{\Phi}}\alpha = 0,$
 i.e., $
\langle P(\alpha),\widetilde{S}_\Omega(\Psi)- \widetilde{S}_\Omega(\Phi)\rangle = 0
$, 
for all closed $1-$form $\alpha$. This implies that $ \widetilde{S}_\Omega(\Psi) 
= \widetilde{S}_\Omega(\Phi)$. 
\end{proof}

Theorem \ref{OrTheo0} implies that if $\Phi$ and $\Psi$ are 
homotopic relatively to fixed endpoints, 
then they have the same flux: Thus, Theorem \ref{OrTheo0}  
implies a result found in \cite{Ban80}. $\maltese$ 
\begin{theorem}\label{OrTheo3} $ $ Let $(M,\Omega)$ be a closed oriented manifold. 
The group $\varGamma_\Omega$ is  
discrete. 
\end{theorem}
\begin{proof} $ $
 The Poincar\'e pairing $\langle,\rangle: H^1(M,\mathbb{R}) 
\times H^{(n-1)}(M,\mathbb{R}) \rightarrow \mathbb{R},$ being 
 a continuous bilinear mapping, let $\mu_0$ denote any   
positive constant such that  
\begin{equation}\label{Disc1}
|\langle A, B \rangle|\leq \mu_0\|A\|_{L^2}\|B\|_{L^2},
\end{equation}
 for 
all $(A, B)\in  H^1(M,\mathbb{R}) 
\times H^{(n-1)}(M,\mathbb{R})$.  Let $\alpha$ (fixed)  be any closed $1-$form such that 
$P(\alpha)\in \overline{\mathcal{B}^1_{\mathcal{S}}(0,1)}$, and consider
 the continuous linear mapping 
$L_\alpha: H^{(n-1)}(M,\mathbb{R}) \rightarrow \mathbb{R}, B\mapsto \langle P(\alpha),  B \rangle$. It 
follows from (\ref{Disc1}) that  $\lvert L_\alpha(B)\rvert  = |\langle P(\alpha), B \rangle| \leq 
\mu_0\lVert B \rVert_{L^2}$, for all 
$B\in H^{(n-1)}(M,\mathbb{R})$ because $ P(\alpha) \in \overline{\mathcal{B}^1_{\mathcal{S}}(0,1)}$. 
Now, since $\varGamma_\Omega$ is countable, we can assume that $\sharp\varGamma_\Omega\textgreater 1$. 
Otherwise, $\varGamma_\Omega =\{ 0\}$, and hence discrete. 
Let $\beta, \tilde\beta\in \varGamma_\Omega$:  By definition 
of $\varGamma_\Omega$, there exist $\Phi, \Psi \in \pi_1(G_\Omega(M))$ such that 
$\beta = \widetilde{S}_\Omega(\Phi),$ and $\tilde\beta = \widetilde{S}_\Omega(\Psi)$. We define another 
loop at the identity by setting $\Theta := \Phi\ast_r\Psi^{-1}$. Therefore, we have following two possibilities. 
 {\bf Either:} There exists a $z_0\in M$ such that 
the $1-$cycle 
$\mathcal{O}_{z_0}^{\Theta} = \mathcal{O}_{z_0}^{\Phi}\bigsqcup(-\mathcal{O}_{z_0}^{\Psi})$ is contractible. In 
this situation, Lemme \ref{Orlem1} implies that
$$\dfrac{1}{Vol_\Omega(M)}\langle P(\alpha), \beta-\tilde\beta \rangle 
= \dfrac{1}{Vol_\Omega(M)}\langle P(\alpha), \widetilde{S}_\Omega(\Theta)\rangle 
= \int_{\mathcal{O}_{z_0}^{\Theta}}\alpha,$$
for all closed $1-$form $\alpha$, i.e.,  
$\beta =\tilde\beta $. {\bf Or:} 
 For each 
$y\in M$, the $1-$cycle $ \mathcal{O}_{y}^{\Theta}$ is not contractible. In this case, since 
 $\beta\neq\tilde\beta$, assume in addition that they are arbitrarily closed 
in the natural topology of $\varGamma_\Omega$, i.e., in the 
topology induced by the vector space structure of $H^{(n-1)}(M, \mathbb{R})$. Since 
 $H^{(n-1)}(M, \mathbb{R})$ is a finite dimensional real vector space, all 
the norms on $H^{(n-1)}(M, \mathbb{R})$ are equivalent, in particular, we shall 
equip $H^{(n-1)}(M, \mathbb{R})$ with the $L^2-$Hodge norm $\lVert, \rVert_{L^2}$: Let 
$l$ be any arbitrary positive integer, and assume that 
\begin{equation}\label{Disc2}
\lVert\beta-\tilde\beta \rVert_{L^2}\textless \frac{ Vol_\Omega(M)}{2\mu_0 l},
\end{equation}
 where $\mu_0$ is the constant defined in (\ref{Disc1}). On the 
other hand,
let $z\in M$ (fixed) and choose any closed $1-$form $\alpha_0$ on $M$ with 
$P(\alpha_0)\in \overline{\mathcal{B}^1_{\mathcal{S}}(0,1)}$ 
such that 
$ \int_{\mathcal{O}_{z}^{\Theta}}\alpha_0 \neq 0$: This is also possible 
because $ \mathcal{O}_{z}^{\Theta}$ is a non-contractible loop by assumption.  
With the help of Lemme \ref{Orlem1} together with (\ref{Disc2}), we derive that
$$0\textless \rvert\int_{\mathcal{O}_{z}^{\Theta}}\alpha_0 \rvert 
= \lvert\dfrac{\langle P(\alpha_0), \widetilde{S}_\Omega(\Theta)\rangle}{Vol_\Omega(M)}\rvert
\leq \frac{\mu_0}{Vol_\Omega(M)}\lVert\beta-\tilde\beta \rVert_{L^2}\leq 1/(2l) \textless 1/l,$$
for all positive integer $l$. This is a contradiction because the top right-hand side 
of the above inequalities tends to zero as $l$ tends to infinity. 
Therefore, we have proved that 
in $\varGamma_\Omega$ two elements are arbitrarily closed, if and only if, they are equal.
\end{proof}

In particular, on any  
Lefschetz closed symplectic manifold $(M,\omega)$, Theorem \ref{OrTheo3} implies that 
the Calabi group is discrete: A result that was definitely  proved by Ono \cite{Ono} for 
any closed symplectic manifold. Furthermore, it is well known that 
the symplectic analogue of Theorem \ref{OrTheo3} which states 
that  the Calabi group is discrete is equivalent to say that  the group of all Hamiltonian diffeomorphisms 
is $C^1-$closed in the group of all symplectic diffeomorphisms isotopic to the identity map. Theorem \ref{OrTheo3} motivated 
the following facts. 
\begin{theorem}\label{RiTheo3} $ $ Let $(M,\Omega)$ be a closed oriented manifold. Then, any 
smooth isotopy $\ker S_\Omega$ is a vanishing-flux path. 
\end{theorem}
\begin{proof} $ $ Let $\Psi =\{\psi^t\}$ be any smooth isotopy in $\ker S_\Omega$, and for each fixed $t$, consider the reparametrized path $\Psi_t:s\mapsto \psi^{(st)}$. Since $M$ is compact, the map $L_\Psi: t \mapsto \Psi_t$ is Lipschitz continuous with respect to the $C^0-$metric, and by Proposition \ref{SC41}, the map $\widetilde S_\Omega$ is continuous with respect to the $C^0-$metric. On the other hand, since $\Psi$ is a  path in $\ker S_\Omega$, then for each $t$,
 we have $\Psi_t(1)\in \ker S_\Omega$, i.e., $\widetilde S_\Omega(\Psi_t)\in \Gamma_\Omega$ for each $t$.  
Thus, the composition $\widetilde S_\Omega\circ L_\Psi$  is a continuous map from the connected topological space $[0,1]
$ into the discrete topological space $\Gamma_\Omega$, hence, this map must be a constant map. That is,  $(\widetilde S_\Omega\circ L_\Psi)(t) = (\widetilde S_\Omega\circ L_\Psi)(0) = 0$, for all $t$. In particular, for $t = 1$, we have\\ $ 0 = (\widetilde S_\Omega\circ L_\Psi)(1) = \widetilde S_\Omega (\Psi)$. 
\end{proof}

\begin{lemma}\label{Rigiditylem0} $ $
 Let $(M,\Omega)$ be a closed oriented manifold. The sub-group $\ker S_\Omega $ 
is $C^1-$closed inside the group $G_\Omega(M)$. 
\end{lemma}

\begin{proof}$ $
 Let $\{\phi_i\}\subset\ker S_\Omega$ 
be a sequence which converges in the $C^1-$metric  
to some $\phi\in G_\Omega(M)$. Let $\{\Phi_i = (\phi_i^t)\}$ be a sequence of 
vanishing-flux isotopies with $\Phi_i\in Iso(\phi_i)_\Omega$ for each $i$, and 
let $\Phi:= (\phi_t)\in Iso(\phi)_\Omega$ be an arbitrary element. Since the 
group $G_\Omega(M) $ is locally connected by smooth arcs (see  \cite{Ban97}), fix a 
 small path connected open neighborhood of the 
identity map, say $\mathcal{V}(id_M)$, and since $\psi^i:=\phi_i\circ\phi^{-1}\xrightarrow{C^1}id_M$, we may 
assume that for all $i$ sufficiently large, one has $\psi^i\in \mathcal{V}(id_M)$. For each $i$ sufficiently 
large, one can choose a smooth curve $\sigma_i$ from $\psi^i$ to the identity, lying inside  $\mathcal{V}(id_M)$ 
such that $\sigma_i$ tends to the constant path identity in the $C^1-$metric as $i$ tends to infinity (if necessary, 
when $\sigma_i$ is only $C^1$, then approximate it in the $C^\infty-$topology by another $\varrho_i$ 
having the same endpoints than $\sigma_i$, \cite{Hirs76}). Assume this done. Now, 
fix any sufficiently large integer  $i$, set $\Psi^i: t\mapsto(\phi_i^t\circ\phi^{-1}_t)$ for each $t$, and  
consider the following juxtaposition of 
paths,
$$
\mathfrak{J}_i(t) =
\left\{
\begin{array}{l}
\Psi^i(\lambda(t)),  \hspace{0,2cm} if \hspace{0,2cm} 0\leq t\leq \frac{1}{2},\\
\sigma_i(\tau(t)),  \hspace{0,2cm} if \hspace{0,2cm}  \frac{1}{2}\leq t\leq 1,    \\
\end{array}
\right.
$$
where $\lambda:[0,1/2]\rightarrow[0,1]$, and $\tau: [1/2,1]\rightarrow [0,1]$ are smooth functions such that 
$\lambda(1/2) = 1 = \tau(1)$, and $\lambda(0) = 0 = \tau(1/2)$. 
By construction, we have $\mathfrak{J}_i(0) = id_M = \mathfrak{J}_i(1)$, i.e., 
$\widetilde{S}_\Omega(\mathfrak{J}_i)\in \varGamma_\Omega$. That is, for all $i$ sufficiently large, 
we have 
$$(\flat): 0= dist(\widetilde{S}_\Omega(\mathfrak{J}_i), \varGamma_\Omega) =
  dist(\widetilde{S}_\Omega(\sigma_i)  - \widetilde{S}_\Omega(\Phi), 
\varGamma_\Omega),$$
 since  $\widetilde{S}_\Omega(\mathfrak{J}_i) = 
\widetilde{S}_\Omega(\sigma_i) + \widetilde{S}_\Omega(\Phi_i) - \widetilde{S}_\Omega(\Phi),$ 
and $\widetilde{S}_\Omega(\Phi_j) = 0$ for all $j$. Taking the 
limit in $(\flat)$ as $i$ tends to infinity yields,  
$$ 0= \lim_i dist(\widetilde{S}_\Omega(\sigma_i)  - \widetilde{S}_\Omega(\Phi), 
\varGamma_\Omega) =  dist(\lim_i\widetilde{S}_\Omega(\sigma_i)  - \widetilde{S}_\Omega(\Phi), 
\varGamma_\Omega)=  dist(- \widetilde{S}_\Omega(\Phi), 
\varGamma_\Omega),$$ because by construction, 
$\sigma_i$ tends to the constant path identity in the $C^1-$metric as $i\rightarrow\infty$. 
That is, $ \widetilde{S}_\Omega(\Phi)\in \overline{\varGamma_\Omega}$, and this implies that 
$\widetilde{S}_\Omega(\Phi)\in \varGamma_\Omega$ because by 
Theorem \ref{OrTheo3}, the group $\varGamma_\Omega$ is  
discrete. Finally, we have proved that any volume-preserving isotopy from the identity 
to $\phi$ must have its flux in  $\varGamma_\Omega$: This means that $\phi\in\ker S_\Omega$. 
\end{proof}

\begin{lemma}\label{Rigiditylem1} $ $
 Let $(M,\Omega)$ be a closed oriented manifold. The sub-group $\ker\widetilde{S}_\Omega $ 
is $C^0-$closed inside the group $Iso(M,\Omega)$. 
\end{lemma}

\begin{proof}$ $

 Let $\{\Phi_i\}$ be a sequence of vanishing-flux isotopies which converges uniformly 
to a volume-preserving isotopy $\Phi$. We have to show that $\widetilde{S}_\Omega(\Phi) = 0$. Since  
 $\Phi$ is smooth, then for each closed $1-$form $\alpha$, we have 
\begin{equation}\label{Rigidityflux3}
 | \int_M \mathcal{F}_{\alpha}(\Phi)(1)\Omega|\leq |\int_M \mathcal{F}_{\alpha}(\Phi_i)(1)\Omega|
+ |\int_M \left(\mathcal{F}_{\alpha}(\Phi_i)(1) -  \mathcal{F}_{\alpha}(\Phi)(1)\right)\Omega| 
\end{equation}
$$\leq 0 +  4Vol_\Omega(M)\|\alpha\|_{L^2} d_{C^0}(\Phi_i,\Phi) \rightarrow 0, i\rightarrow\infty,$$
because by Proposition \ref{SC41} we have,  
$0 = \langle P(\alpha), \widetilde{S}_\Omega(\Phi_i)\rangle = \int_M \mathcal{F}_{\alpha}(\Phi_i)(1)\Omega,$ 
for each $i$, and since $\lim_{C^0}(\Phi_i) =\Phi$, it follows from Lemma $3.10-$\cite{TD2} that 
$$ |\int_M \mathcal{F}_{\alpha}\left((\Phi_i)(1) -  \mathcal{F}_{\alpha}(\Phi)(1)\right)\Omega| \leq 
Vol_\Omega(M)\sup_{x\in M}|\mathcal{F}_{\alpha}(\Phi_i)(1)(x) -\mathcal{F}_{\alpha}(\Phi)(1)(x)|,$$
$$ \leq 4Vol_\Omega(M)\|\alpha\|_{L^2}d_{C^0}(\Phi_i,\Phi),$$
for $i$ large enough. Thus, the right-hand side in (\ref{Rigidityflux3}) tends to zero as $i$ goes 
at infinity, and hence the left-hand side in (\ref{Rigidityflux3}) vanishes. Therefore, since 
$\Phi$ is a volume-preserving isotopy, it follows from Proposition \ref{SC41} that, 
$ 0 = \int_M \mathcal{F}_{\alpha}(\Phi)(1)\Omega = \langle P(\alpha), \widetilde{S}_\Omega(\Phi)\rangle,$ 
for each closed $1-$form $\alpha$. That is, 
$ \langle P(\beta), \widetilde{S}_\Omega(\Phi)\rangle = 0,$ for all closed $1-$form $\beta$. This, 
implies that $\widetilde{S}_\Omega(\Phi) = 0$.
\end{proof} 
\subsection{Application to quasi-morphisms for symplectic closed $2-$surfaces}\label{Quasi-Mor} 
We assume that $(M, \omega)$ is a $2-$dimensional closed symplectic manifold. Firstly, 
observe that Proposition \ref{GF1} implies that for each non-trivial closed $1-$form 
$\alpha,$ and each $x\in M$, we have
\begin{equation}\label{QS1}
 \varDelta(\psi\circ\phi, \alpha)_x = \varDelta(\psi, \alpha)_x 
+ \varDelta(\phi, \alpha)_x 
+ \dfrac{Vol_\Omega(M)}{\|\alpha\|_{L^2}}\left(\int_{\mathcal{O}_{x}^\Psi}\alpha- 
\int_{\mathcal{O}_{\phi(x)}^\Psi}\alpha\right),
\end{equation}
for all $\psi, \phi \in Symp_0(M,\omega)$, and for all $\Psi\in Iso(\psi)_\omega$.  
That is, the quantity\\  $\left(\int_{\mathcal{O}_{x}^\Psi}\alpha- 
\int_{\mathcal{O}_{\phi(x)}^\Psi}\alpha\right)$ 
seems to impeach $\varDelta(., \alpha)_x$ to be 
a homomorphism in the usual sense. However, the goal of this 
subsection is to check whether $\varDelta(., \alpha)_x$ can give rise to a 
 quasi-morphism or not.\\ Note that if $G$ is a group, then 
 a mapping $\varDelta: G \rightarrow \mathbb{R}$ is a 
quasi-morphism if there exists a real number $R\textgreater 0$ such that
$|\varDelta(fg) -\varDelta(f) -\varDelta(g)|\leq R,$
for all $f,g \in G$. In addition, if $\varDelta(f^m) = m\varDelta(f),$ 
for each $f\in G$ an for each $m\in \mathbb{Z}$, then 
 $\varDelta$ is called a homogeneous quasi-morphism (see \cite{CB}). 
The minimal value of the real number $R$ in the above definition 
is sometimes called the defect of the quasi-morphism $\varDelta$. This 
 is due to the fact that if $R = 0$, then  $\varDelta$ will be a homomorphism in the usual sense. 
Furthermore, for each quasi-morphism $\varDelta$, the map $\psi 
\mapsto \lim_{m\longrightarrow\infty} \{\varDelta(\psi^m)/m\} ,$ induces a homogeneous quasi-morphism.
\begin{proposition}\label{Qusi-pro} $ $
Assume $(\widetilde M, \omega)$ to be a $2-$dimensional closed 
symplectic manifold with $H^1(\widetilde M,\mathbb{R}) \neq 0$.
Let $\alpha$ be a non-trivial closed $1-$form and $p\in \widetilde M$ be any point. Therefore, the mapping 
 $\varDelta( . ,\alpha)_p : 
Symp_0(\widetilde M,\omega) \rightarrow \mathbb{R}, \psi \mapsto \varDelta(\psi, \alpha)_p,$
is a continuous quasi-morphism whose defect is less than or equal 
to twice the square of the symplectic area of $\widetilde M$.
\end{proposition}
\begin{proof} $ $
 
Let  $\alpha$ be a non-trivial closed $1-$form and 
 $p\in \widetilde M$. Proposition \ref{GF1} implies that for all $\psi, \phi \in Symp_0(\widetilde M,\omega)$, we have 
$$\dfrac{1}{Vol_\Omega(\widetilde M)}|\varDelta(\psi\circ\phi, \alpha)_p 
- \varDelta(\psi, \alpha)_p - \varDelta(\phi, \alpha)_p |
  \leq |\int_{\mathcal{O}_{\phi(p)}^\Psi}\dfrac{1}{\|\alpha\|_{L^2}}\alpha| 
+ |\int_{\mathcal{O}_{p}^\Psi}\dfrac{1}{\|\alpha\|_{L^2}}\alpha |,$$
for each $\Phi\in Iso(\phi)_\omega$. Since $\int_{\mathcal{O}_{p}^\Psi}\dfrac{1}{\|\alpha\|_{L^2}}\alpha $ 
(resp. $\int_{\mathcal{O}_{\phi(p)}^\Psi}\dfrac{1}{\|\alpha\|_{L^2}}\alpha$) 
is the algebraic value of the symplectic area swept 
by the orbit $\mathcal{O}_{p}^\Psi $ (resp. $\mathcal{O}_{\phi(p)}^\Psi $) 
under the symplectic flow generated by the 
closed $1-$form 
$\dfrac{1}{\|\alpha\|_{L^2}}\alpha$, then 
$|\int_{\mathcal{O}_{p}^\Psi}\dfrac{1}{\|\alpha\|_{L^2}}\alpha|$ 
(resp. $|\int_{\mathcal{O}_{\phi(p)}^\Psi}\dfrac{1}{\|\alpha\|_{L^2}}\alpha|$) is less than or equal 
to the symplectic area $A(\widetilde M)$ of $\widetilde M$. Therefore, 
\begin{equation}
 |\varDelta(\psi\circ\phi, \alpha)_p - \varDelta(\psi, \alpha)_p 
- \varDelta(\phi, \alpha)_p | \leq 2 A(\widetilde M)Vol_\Omega(\widetilde M) = 2[A(\widetilde M)]^2,
\end{equation}
for all $\psi, \phi \in Symp_0(\widetilde M,\omega)$ because here the area form is the volume form. 
\end{proof}

\section{On the dynamics of fix-points for volume-preserving diffeomorphisms}\label{Fixpoints}
Let $\alpha \in \mathcal{Z}^1(M)$, and fix a point $x\in M$.  
 Corollary \ref{GF10} implies that for each $\varphi\in G_\Omega(M),$
the real number $\varDelta( \varphi ,\alpha)_x$ does not depend on any choice of an isotopy  
$\Psi$ in $Iso(\varphi)$. In particular, as we shall see this seems to suggest that the 'zeros' of the function 
 $\varDelta(.,\alpha)_x$ could fix-points for volume-preserving diffeomorphisms.
The goal of this section is to enhance the application  
of flux geometry in the study of the dynamics of fix-points for volume-preserving diffeomorphisms. 
To that end, let us introduce the following notions and notations.\\

 Given a non-empty set $\mathfrak{A}$, we shall write  $\sharp \mathfrak{A}$ 
to mean the cardinal of $\mathfrak{A}$; while
 for each $\phi\in G_\Omega(M)$, we shall write $\mathcal{FIX}(\phi)$ to mean the set of 
all point $z\in M$ such that $\phi(x) = x$.\\ 
Now, let $h\in C^\infty(M,\mathbb{R})$ be non-trivial, and let $Crit(h)\subset M$, be the subset of $M$ 
consists of all the critical points for $h$. From Morse theory, we may assume that $Crit(h)\neq \emptyset$,  
 for each $h\in C^\infty(M,\mathbb{R})$. In fact, the Morse 
theory seems to suggest that  critical points of smooth functions on a closed 
oriented manifold are potential candidates to be fix-points for volume-preserving diffeomorphisms 
in the following sense.\\ 
{\bf Claim:} If one assumes that there is $\psi \in G_\Omega(M)$ such that $\mathcal{FIX}(\psi) =\emptyset$, then 
no smooth function on $M$ can be injective.\\ 
{\bf Proof of the claim:} If $h\in C^\infty(M,\mathbb{R})$ was injective,
 then Proposition \ref{SC41} will imply that for each $\Phi\in Iso(\psi)$, we could have  
$
 \int_M\mathcal{F}_{dh}(\Phi)\Omega = \langle  P(dh), \widetilde{S}_\Omega(\bar{\Phi})\rangle = 0,$
because $ P(dh) = 0$. Hence, we could have a point $x(h,\Phi)\in M$  (depending on 
 $h$ and $\Phi$) such that $\mathcal{F}_{dh}(\Phi)(x(h,\Phi)) = 0$, i.e.,
\begin{equation}\label{ExpFix}
 0 = \mathcal{F}_{dh}(\Phi)(x(h,\Phi)) =\int_{\mathcal{O}_{x(h,\Phi)}^\Phi}dh = h(\psi(x(h,\Phi))) - h(x(h,\Phi)),
\end{equation}
and the injectivity of $h$ gives, $\psi(x(h,\Phi))) = x(h,\Phi)$. So, the injective nature of $h$ would be in  
contradiction the fact that $\mathcal{FIX}(\psi) =\emptyset$. $\maltese$ 
\begin{example}

Consider the torus $\mathbb T^{2k}$ with coordinates $(\theta_1,\dots,\theta_{2k})$ and 
equip it with the flat Riemannian metric $g_0$:  
All the $1-$forms $d\theta_i$, $i= 1,\dots, 2k$ are harmonic. 
Take the $1-$forms $d\theta_i$ for $i= 1,\dots, 2k$ as basis for the space of harmonic $1-$forms 
and consider the symplectic form
$
\omega = \sum_{i=1}^k d\theta_i \wedge d\theta_{i +k}.
$ 
 Then, the action 
of the unit circle $\mathbb S^{1}$ on $\mathbb T^{2k}$:\\ 
$ \rho: \mathbb S^{1}\times\mathbb T^{2k}\rightarrow \mathbb T^{2k}, 
(\alpha, (\theta_1, \theta_2, \dots,\theta_{2k}))\mapsto 
(\theta_1 + \alpha, \theta_2 + \alpha, \dots,\theta_{2k} + \alpha),$ 
induces a volume-preserving diffeomorphism $\rho_\alpha: \mathbb T^{2k}\rightarrow \mathbb T^{2k}$ 
which has no fixed point whenever $\alpha$ is small and non-trivial. 
Thus, no smooth function in $C^\infty(\mathbb T^{2k},\mathbb{R})$ can be injective. 
\end{example}



An interesting problem in fix-points theory is how to predict exactly the minimum 
number of fix-points that can have a mapping. In the 
frame of symplectic geometry, 
we have the following conjecture which is due to Arnold.\\

{\bf Conjecture (Arnold):}
A Hamiltonian diffeomorphism of a closed symplectic manifold $(M,\omega)$ must have at least 
as many fixed points as the minimal number of critical points of a smooth function on $M$.\\

Here are some affirmatives answers to the Arnold conjecture:
\begin{itemize}
 \item Using Lagrangian intersections, Weinstein \cite{Wien1} showed that: 
Any symplectic diffeomorphism
of a  compact simply connected symplectic manifold has at least two fixed-points
provided that it is sufficiently $C^1-$closed to the identity map.
\item Banyaga \cite{Ban80}: Showed that 
any Hamiltonian diffeomorphism has at least $cat(M)$ fixed-points
provided that it is sufficiently $C^1-$closed to the identity map, 
where $cat(M)$ stands for the Lusternik-Schnirelman category 
of the manifold $M$ \cite{Ban80}: $cat(M)$ is the smallest positive integer $l$ such that 
there exists a cover $U_1,\dots,U_{l + 1}$ of $M$ consists of contractible open subsets.
 Lusternik-Schnirelman proved that the inequality holds\\ 
$1+ cat(M)\leq \sharp Crit(f)\leq 1 + \dim(M)$, 
for any smooth function $f$ on $M$.
\item Definitely, the Arnold conjecture was solved by Fukaya-Ono 
\cite{F-O}, and Liu-Tian \cite{L-T}, using different methods.
\end{itemize}

We have the following consequence of Proposition \ref{SC41}.

\begin{proposition}\label{FiXPRO} $ $
 Let $(M,\Omega)$ be an oriented closed manifold, and let $\psi\in G_\Omega(M)$. Assume that 
there exists $h\in C^\infty(M,\mathbb{R})$ such that $\sharp Crit(h)= 1 $. If there exists 
 $\Phi\in Iso(\psi)_\Omega$ such that $x(h,\Phi)\in Crit(h)$ 
where the point $x(h,\Phi)$ is defined as in (\ref{ExpFix}), then $x(h,\Phi)$ is a fixed-point for $\psi$.
\end{proposition}
\begin{proof} $ $
 Let $\psi\in G_\Omega(M)$, assume that 
there exists $h\in C^\infty(M,\mathbb{R})$ such that $\sharp Crit(h)= 1$, and 
that there exists $\Phi\in Iso(\psi)_\Omega$ such that $x(h,\Phi)\in Crit(h)$.  
Set $z:= x(h,\Phi)$, and assume that $\psi(z) \neq z$. By assumption, since $\psi(z) \neq z$, then 
the  smooth map 
$t\mapsto h(\Phi(z,t))$ satisfies exactly one of the following conditions: Either, 
\begin{equation}\label{IN1}
 h(z)\textless h(\psi(z)),
\end{equation}
 or 
\begin{equation}\label{IN2}
 h(\psi(z))\textless h(z).
\end{equation}
 On the other hand, since $\int_{\mathcal{O}_z^\Phi}dh = 0$ (Proposition \ref{SC41}), 
then we derive that $ h(\psi(z))= h(z)$. This is in contradiction with any of the inequalities (\ref{IN1}) and (\ref{IN2}).
 Thus, 
$\psi(z) = z$. 
\end{proof}
\begin{definition}
An element $(p,\psi)\in M\times G_\Omega(M)$ is called  a zero for $\varDelta$, if and only 
if, $\varDelta(\psi, \alpha)_p = 0$, for all $\alpha \in \mathcal{Z}^1(M)$. 
\end{definition}

\begin{theorem}\label{Functor3} $ $
 Let $(M,\Omega)$ be an oriented closed manifold. Then, $ (p,\psi)\in M\times\ker S_\Omega$ is a
zero for $\varDelta$, if and only 
if, $p\in \mathcal{FIX}(\psi)$.
\end{theorem}
The proof of Theorem \ref{Functor3} will need the following contractibility result which  generalizes a result 
from Hamiltonian dynamics found in \cite{McDuff-SAl, TD2}. 
\begin{lemma}\label{FixPt4}{\bf (Global contractibility)} $ $
 Let $(M,\Omega)$ be an oriented closed manifold. If $\psi\in \ker S_\Omega$ 
has a fix-point $p$, then for any $\Phi\in Iso(\psi)_\Omega$ with trivial flux,  
the orbit $ \mathcal{O}_p^\Phi$ is contractible.
\end{lemma}

Here is a simple case of Lemma \ref{FixPt4}.
\begin{corollary}\label{0FixPt4}{\bf (Local contractibility)}
 Consider the set 
$$C_{r(g)}(id_M): =\{\phi\in G_\Omega(M): d_{C^1}(\phi,id_M)\textless r(g)\},$$ where 
$r(g)$ is the injectivity radius 
of a Riemannian metric $g$ on $M$. Then, for each $p\in \mathcal{FIX}(\psi) $, and for any 
isotopy $\Phi\in Iso(\psi)_\Omega$ which does not escape $C_{r(g)}(id_M)$, the orbit $ \mathcal{O}_p^\Phi$ 
is contractible. 
\end{corollary}
\begin{proof} $ $
 Let $\Phi=\{\phi_t\}\in Iso(\psi)_\Omega$ which does not escape $C_{r(g)}(id_M)$. Then, the condition 
$ d_{C^1}(\phi_t,id_M)\textless r(g)$, for all $t$, implies in particular that 
for each $t$, the points $\phi_t(p)$ can be connected to $p$ via a unique minimal geodesic $C_{\phi_t(p),p}$. 
This induces a homotopy between the curves $ \mathcal{O}_p^\Phi$  and $t\mapsto p$. We may define this homotopy $H$
as follows: $H(0,0) = p = \phi(p) = H(1,1)$, $H(0,t) = \mathcal{O}_p^\Phi$, $H(1,t) = p$, for all $t$, and 
$H(s,t)$ is the time$-s$ of the geodesic $C_{\phi_t(p),p}$, i.e., 
$\circledast(\Phi, p):=\{H(s,t): (s,t)\in [0,1]\times [0,1]\}$ is a $2-$chain with 
$\partial\circledast(\Phi, p) = \mathcal{O}_p^\Phi\sqcup\{p\}= \mathcal{O}_p^\Phi$.
\end{proof}

\begin{corollary}\label{1FixPt4} $ $
 Let $\{\mathbf{U}_i\}_{1\leq i \leq k}$ be an open cover of $M$. Let $\psi\in \ker S_\Omega$ 
which can be fragmented as $\psi = \psi_1\circ\dots\circ\psi_k$ such that for each $i$:
\begin{itemize}
 \item  $\psi_i\in \ker S_\Omega$, and $supp(\psi_i)\subset \mathbf{U}_i $,  
\item $\psi_i$ is sufficiently $C^1-$closed to the identity map,
\item $\psi_i$ can be connected to the 
identity map through a vanishing-flux path $\Phi_i$ which is sufficiently $C^1-$closed to the identity 
path.
\end{itemize}
Then,  for each  $\alpha \in \mathcal{Z}^1(M)$, we have $\varDelta(\psi, \alpha)_p = 0$, 
provided $p\in \mathcal{FIX}(\psi_i) $.
\end{corollary}
\begin{proof}$ $
   Since $\psi = \psi_1\circ\dots\circ\psi_k$, one derives from  Proposition \ref{GF1} 
that 
$$\varDelta(\psi, \alpha)_p = \varDelta(\psi_1\circ\dots\circ\psi_k, \alpha)_p = 
\varDelta(\psi_k, \alpha)_p + 
\sum_{k-1\geq i\geq 1}\varDelta(\psi_i, \alpha)_{(\psi_{i+ 1}\circ\dots\circ\psi_k)(p)},$$
for each  $\alpha \in \mathcal{Z}^1(M)$. On the other hand, since each $\psi_i\in \ker S_\Omega$ is sufficiently 
$C^1-$closed to the identity map, and $\ker S_\Omega$ is locally 
connected by smooth arcs, let $\Phi_i\in Iso(\psi_i)_\Omega$ be an isotopy with trivial flux 
which is   
sufficiently $C^1-$closed to the identity path. For each non-trivial $\alpha \in \mathcal{Z}^1(M)$, 
Corollary \ref{GF10} yields: $\|\alpha\|_{L^2}\varDelta(\psi_k, \alpha)_p = 0 
-Vol_\Omega(M)\int_{\mathcal{O}^{\Phi_k}_p}\alpha = 0,$ 
because $p$ being a fix-point for $\psi_k$, and the vanishing-flux isotopy $\Phi_k$ being  
sufficiently $C^1-$closed to the identity path, it follows from Corollary \ref{0FixPt4}, that 
the orbit $\mathcal{O}^{\Phi_k}_p$ 
is contractible.  For each $i \textless k$, we have\\ 
$\|\alpha\|_{L^2}\varDelta(\psi_{i}, \alpha)_{(\psi_{i+ 1}\circ\dots\circ\psi_k)(p)} = 
\|\alpha\|_{L^2}\varDelta(\psi_{i}, \alpha)_p 
= 0 - Vol_\Omega(M)\int_{\mathcal{O}^{\Phi_{i}}_p}\alpha = 0,$
because $p$ being a fix-point for $\psi_{i}$, and the vanishing-flux isotopy $\Phi_{i}$ being 
sufficiently $C^1-$closed to the identity path, it follows from Corollary \ref{0FixPt4}, that 
the orbit $\mathcal{O}^{\Phi_{i}}_p$ 
is contractible. 
That is, we have $\varDelta(\psi, \alpha)_p = 0,$ for all $\alpha \in \mathcal{Z}^1(M)$.
\end{proof}

\begin{remark}\label{2FixPt4} 
\begin{itemize}
 \item The hypotheses of Corollary \ref{1FixPt4} are all fulfilled when $\Omega$ is the 
symplectic  volume form.
\item If we let $\psi$ as in Corollary \ref{1FixPt4}, then by the help of Corollary \ref{GF10} we derive that 
for any fix-point $p$ 
for $\psi$ (whenever it exists) and 
for any vanishing-flux path $\Phi\in Iso(\psi)_\Omega$, the orbit $\mathcal{O}^{\Phi}_p$ 
is contractible. Furthermore, in this context with the help of Corollary \ref{GF10}, the $1-$cycle 
$\mathcal{O}^{\Phi}_p$ is homologous to a $1-$cycle obtained by gluing at $p$, the 
contractible $1-$cycles $\mathcal{O}^{\Phi_i}_p$: We shall call the latter $k-$fold  of 
$1-$cycles, 
a $k-$jet of contractible $1-$cycles at $p$.
\end{itemize}

\end{remark}

{\it Proof of Lemma \ref{FixPt4} ~:} Let $\phi\in \ker S_\Omega$  which 
has a fix-point $p$. 
\begin{itemize}
	\item {\large Step$(1/3)$}: 
\end{itemize}
Note that by Lemma \ref{Orlem1} if the orbit of a fix point $p$ of $\phi$ under a particular vanishing-flux Hamiltonian path with time-one map $\phi$ is contractible, then so is the orbit of $p$ under any other vanishing-flux path with time-one map $\phi$. 
So, we shall proceed by contradiction by assuming that for 
any vanishing-flux path $\Phi$ from the identity to $\phi$, the orbit $\mathcal{O}_p^{\Phi}$ is not contractible. 
\begin{itemize}
	\item {\large Step$(2/3)$}: 
\end{itemize}
 Consider any vanishing-flux path $t\longmapsto \Phi_t$ with 
time-one map $\phi$, and let $\psi\in \ker S_\Omega \backslash \{\phi\}$ be another vanishing-flux map. Since 
Theorem \ref{RiTheo3} shows that any  isotopy  in $\ker S_\Omega$ is a vanishing-flux path, while  the group $\ker S_\Omega$ is paths connected through smooth paths, then 
 consider $t\longmapsto H_t$ to be a vanishing flux path from $\phi$ to the identity, and $t\longmapsto \Xi_t$ to be any vanishing-flux path from $\phi$ to $\psi$. Consider the following juxtaposed paths: 
$$
\zeta_t  =
\left\{
\begin{array}{l}
 H_{(2t)},  \hspace{0,2cm} if \hspace{0,2cm} 0\leq t\leq \frac{1}{2},\\
\Xi_{ (2t-1)},  \hspace{0,2cm} if \hspace{0,2cm}  \frac{1}{2}\leq t\leq 1, \\
\end{array}
\right.
$$\\
and,
$$
L_t  =
\left\{
\begin{array}{l}
\Phi_{(2t)},  \hspace{0,2cm} if \hspace{0,2cm} 0\leq t\leq \frac{1}{2},\\
\zeta_{(2t -1)},  \hspace{0,2cm} if \hspace{0,2cm}\frac{1}{2} \leq t\leq 1.\\
\end{array}
\right.
$$\\
 Since $L$ is a vanishing-flux loop at the identity map (if necessary one can smooth $L$ where there is need), then 
Lemma \ref{Orlem1}, all the orbits of $L$ are contractible, in particular, the loop 
$\mathcal{O}_p^{L},$ is contractible, i.e., the loop $\mathcal{O}_p^{\Phi}\cup\mathcal{O}_{p}^{\zeta},$ is contractible. 
\begin{itemize}
	\item {\large Step$(3/3)$}: 
\end{itemize}
On the other hand, since the parametrized path $\bar\gamma: t\longmapsto \zeta_{(1-t)}$, is a vanishing-flux path with time-one map $\phi$, then by assumption, the orbit $\mathcal{O}_p^{\bar\gamma}$ is not contractible, and so is the orbit $\mathcal{O}_p^{\zeta}$.  Thus, the orbit $\mathcal{O}_p^{L},$ is 
made up by two non-contractible loops $ \mathcal{O}_p^{\Phi}$, and $\mathcal{O}_p^{\zeta}$. Since both loops $ \mathcal{O}_p^{\Phi}$, and $\mathcal{O}_p^{\zeta}$ have the 
same orientation, then the loop $  \mathcal{O}_p^{\Phi}\cup\mathcal{O}_p^{\zeta}$ cannot be contractible: 
a particle moving along $\mathcal{O}_p^{L},$ will have to travel from $t = 0$ to $t =1/2$ along the loop $\mathcal{O}_p^{\Phi},$ and then from $t = 1/2$ to $t = 1$, the same particle 
will be continue to move along $\mathcal{O}_p^{\zeta}$ without reversing the orientation $\mathcal{O}_p^{\Phi}$.  Therefore,  {\large Step$(2/3)$}  and {\large Step$(3/3)$} contradict themself. Q.E.D\\

{\it Proof of Theorem \ref{Functor3}~:} Assume that
$ (p,\psi)\in M\times\ker S_\Omega$ is a
zero of $\varDelta$. Pick $\Phi\in \left(Iso(\psi)_\Omega\cap\ker\widetilde{S}_\Omega\right)$, and 
derive with the help of Corollary \ref{GF10} that\\ 
$ 0= \|\alpha\|_{L^2}\varDelta(\psi, \alpha)_p =  0 - Vol_\Omega(M)\int_{\mathcal{O}_p^\Phi}\alpha$, 
for each $\alpha \in \mathcal{Z}^1(M)\backslash\{0\}$, i.e., $ \int_{\mathcal{O}_p^\Phi}\alpha  = 0,$ 
for all closed $1-$form $\alpha$. Thus, the orbit $\mathcal{O}_p^\Phi $ is a boundary, i.e.,
 $\psi(p) = p$. Conversely, 
assume that $\psi(p) = p$: Since $\psi\in \ker S_\Omega$, then 
there is an isotopy $\Phi\in Iso(\psi)_\Omega$ such that $\widetilde{S}_\Omega(\Phi) = 0.$ So, 
 Corollary \ref{GF10} implies that 
 $\|\alpha\|_{L^2}\varDelta(\psi, \alpha)_p =  - Vol_\Omega(M)\int_{\mathcal{O}_p^\Phi}\alpha$
 for each $\alpha \in \mathcal{Z}^1(M)\backslash\{0\}$, i.e., 
 $\varDelta(\psi, \alpha)_p = 0$, 
for each $\alpha \in \mathcal{Z}^1(M)$, because 
$p$ being a fix-point for $\psi$, it follows from
 Lemma \ref{FixPt4} that the orbit $\mathcal{O}_p^\Phi $ is contractible. $\square$\\

The following fix-points result is 
an attempt to a generalization of the 
Arnold conjecture: A consequence of the Thurston fragmentation result together with Proposition \ref{GF1}. 
 
\begin{lemma}\label{Functor7}{\bf (Weak Arnold's conjecture)}  $ $Let $(M,\Omega)$ be an oriented closed manifold. 
Any vanishing-flux volume-preserving diffeomorphism which is sufficiently $C^1-$closed 
to the identity must have at least 
as many fixed points as the minimal number of critical points of a smooth function on $M$.
\end{lemma}
\begin{proof} $ $
 
 Since $M$ is compact, assume that 
 $\psi = \psi_1\circ\dots\circ\psi_k$ with respect to some open cover $\{\mathbf{W}_i\}_{1\leq i\leq k}$ of $M$ 
with 
each $\psi_i\in \ker S_\Omega$ supported in $\mathbf{W}_i$. We proceed by contradiction: 
Assume that $\psi(x)\neq x$, for all $x\in M$. This means that any $x\in M$ belongs to the support of a certain 
$\psi_i$, $i=1,\dots,k$, i.e. $\psi_i(x)\neq x$, for all $x\in M$, and for each $i=1,\dots,k$.  
This fragmentation 
together with Proposition \ref{GF1} implies that,
\begin{equation}\label{0Arnold}
 \varDelta(\psi, \alpha)_p = \varDelta(\psi_1\circ\dots\circ\psi_k, \alpha)_p = 
\varDelta(\psi_k, \alpha)_p + 
\sum_{k-1\geq i\geq 1}\varDelta(\psi_i, \alpha)_{(\psi_{i+ 1}\circ\dots\circ\psi_k)(p)},
\end{equation}
for each  $\alpha \in \mathcal{Z}^1(M)$, and all $p\in M$. Since by assumption 
$\psi$ is 
 sufficiently $C^1-$closed 
to the identity map, we can also assume that there exists an integer $i_0$ such that $\psi_{i_0}$ is arbitrary 
$C^1-$closed 
to the identity map. For such an integer $i_0$, pick $y_0\in supp(\psi_{i_0})$ with $y_0\notin supp(\psi_{j})$ 
for all $j\neq i_0$, and since we must have 
 $y_0 \neq\psi(y_0)$, then construct a smooth function $h\in C^\infty(M,\mathbb{R})$ with the following properties:
 $ supp(h) \varsubsetneq supp(\psi_{i_0})$,  
 $y_0\in supp(h)$,  $\psi_{i_0}(y_0)\notin supp(h)$,
 and  $h = 3$,
on a small open 
neighborhood $\mathcal{K}(y_0)\subset supp(h)$ of $y_0$. 
As in the proof of Lemma \ref{Cont1}, since $\psi_{i_0}$ is arbitrary 
$C^1-$closed 
to the identity map we may have the following inequality,\\  $(\ast): 
  |\varDelta(\psi_{i_0}, dh)_{(\psi_{i_0 + 1}\circ\dots\circ\psi_k)(y_0)}| 
 = |\varDelta(\psi_{i_0}, dh)_{(y_0)}| \leq 2 Vol_\Omega(M)d_{C^1}(id_M,\psi_{i_0}).$ 
 On the other hand,  we also compute 
 $ \varDelta(\psi_j, dh)_{(\psi_{j + 1}\circ\dots\circ\psi_k)(y_0)} = h(y_0) 
-  h(y_0) = 0,$
whenever $j \textgreater i_0$, and 
  $ \varDelta(\psi_j, dh)_{(\psi_{j + 1}\circ\dots\circ\psi_k)(y_0)} = h(\psi_{i_0}(y_0)) -
 h(\psi_{i_0}(y_0)) = 0,$
whenever $j\textless i_0$. 
 Clearly, (\ref{0Arnold}) implies that $
  \varDelta(\psi, dh)_{y_0} = \varDelta(\psi_{i_0}, dh)_{y_0}$, and 
combining this last equality together with ($\ast$) yields,
$ 1 \textless \frac{1}{2}|h(y_0)|\leq \|dh\|_{L^2}d_{C^1}(id_M,\psi_{i_0})$:   
This is not plausible because the top right-hand side can be 
considered as arbitrarily small, while the top left-hand 
side remains constant.  Of course, it is not hard to see that $y_0$ is
 a critical point for $h$.
\end{proof}

Lemma \ref{Functor7} implies the following fix-points result which generalizes and solves the 
Arnold conjecture. 

\begin{theorem}\label{Functor11}{\bf (Strong Arnold's conjecture)} $ $
 Let $(M,\Omega)$ be an oriented closed manifold. 
Any vanishing-flux volume-preserving diffeomorphism must have at least 
as many fixed points as the minimal number of critical points of a smooth function on $M$.
\end{theorem}
\begin{proof} $ $
 If $\psi\in \ker S_\Omega$ is sufficiently small in the $C^1-$metric, then apply Lemma \ref{Functor7}
to conclude. Otherwise, $\psi$ can be fragmented as  
 $\psi_1\circ\dots\circ\psi_k$ with respect to some open cover $\{\mathbf{V}_i\}_{1\leq i\leq k}$ of $M$ with 
 $\psi_i\in \ker S_\Omega$ and $supp(\psi_i)\subset\mathbf{V}_i$. 
If necessary, repeat this fragmentation process until 
to obtain a term $\psi_j$ which is sufficiently $C^1-$closed to the identity map, and conclude 
by applying Lemma \ref{Functor7}. 
\end{proof}

The following result shows that a fix-point of any $\psi\in\ker S_\Omega$  
can determine the flux geometry of any path in $Iso(\psi)_\Omega$. 

\begin{lemma}\label{Functor55}
 Let $(M,\Omega)$ be an oriented closed manifold. Let $\phi\in \ker S_\Omega$. If $\Psi\in Iso(\phi)_\Omega$, then 
we must have 
 $$\langle P(\alpha), \widetilde{S}_\Omega(\Psi)\rangle = Vol_\Omega(M)\int_{\mathcal{O}_p^\Psi}\alpha,$$
 for each $\alpha \in \mathcal{Z}^1(M)$, whenever $p$ is a fix point for $\phi$. 
\end{lemma}
\begin{proof} $ $
 
 Let $\Psi\in Iso(\phi)_\Omega$, pick $\Phi\in Iso(\phi)_\Omega$ with vanishing-flux,
 and derive from Corollary \ref{GF10} that 
$$\langle P(\alpha), \widetilde{S}_\Omega(\Psi)\rangle
  -Vol_\Omega(M)\int_{\mathcal{O}_p^\Psi}\alpha
= \|\alpha\|_{L^2}\varDelta(\phi, \alpha)_p = \langle P(\alpha), \widetilde{S}_\Omega(\Phi)\rangle
 -Vol_\Omega(M)\int_{\mathcal{O}_p^\Phi}\alpha,$$
for each non-trivial $\alpha \in \mathcal{Z}^1(M)$, i.e., $
 \langle P(\alpha), \widetilde{S}_\Omega(\Psi)\rangle
  -Vol_\Omega(M)\int_{\mathcal{O}_p^\Psi}\alpha
= 0,$ 
for each $\alpha \in \mathcal{Z}^1(M)$, because $ \widetilde{S}_\Omega(\Phi) = 0$, and by 
  Lemma \ref{FixPt4} the orbit $\mathcal{O}_z^\Phi$ is contractible whenever $z\in  \mathcal{FIX}(\phi)$. 
\end{proof}

\subsection{Fix-points and $C^0-$dynamics}
The techniques used in the proof of the above fix-points results require the mapping to be at least 
of class $C^1$. In what following, we wish to experiment another way for searching 
fix-points for volume-preserving diffeomorphisms using only 
the $C^0-$metric. For this purpose we will need the following notions:\\ 
For each $\alpha \in \mathcal{Z}^1(M)$, and each $\psi\in G_\Omega(M)$, 
let $\mathcal{N}^\varDelta(\psi,\alpha)\subset M$ denote the null set of the 
function $p\mapsto\varDelta(\psi, \alpha)_p$. 
It is not hard to see that  for each 
$(\psi, \alpha) \in G_\Omega(M)\times \mathcal{Z}^1(M)$, the set
$ \mathcal{N}^\varDelta(\psi,\alpha)$ is closed and non-empty.

\begin{corollary}\label{Fix-Points3-5} $ $
 Let $(M,\Omega)$ be a closed symplectic manifold. If $\psi\in\ker S_\Omega$ is sufficiently  
$C^0-$closed to the constant map identity,  
then 
$\bigcap_{\alpha\in \mathcal{Z}^1(M)}\mathcal{N}^\varDelta(\psi,\alpha)\neq\emptyset.$
\end{corollary}
\begin{proof}
 
 Let $\psi\in G_\Omega(M)$ sufficiently 
$C^0-$closed to the constant map identity. Let $\alpha$ 
be an arbitrary closed $1-$form. Since by Lemma \ref{Cont1}, we have\\ 
$$\varDelta(\psi, \alpha)\xrightarrow{\psi\xrightarrow{C^0}id_M} \varDelta(id_M, \alpha) 
= 0\in C^\infty(M,\mathbb{R}),$$
 then it follows that the closed subset 
$\mathcal{N}^\varDelta(\psi,\alpha)$ tends to $M$ whenever $\psi\xrightarrow{C^0}Id_M$ (no matter 
the choice of $\alpha$). 
So, for each finite family $\{\alpha_i\}_{i\in I}$ of closed $1-$forms, the set 
$ \bigcap_{ i\in I}\mathcal{N}^\varDelta(\psi,\alpha_i)$ can be assumed to be non-empty
 whenever $\psi$ is assumed to be sufficiently 
$C^0-$closed to the identity map because under this assumption, each  
closed subset $\mathcal{N}^\varDelta(\psi,\alpha_i)$ 
tends to recover the ambient manifold $M$.  Hence, 
 the  family 
of closed subsets $\{\mathcal{N}^\varDelta(\psi,\alpha)\subset M: 
\alpha\in \mathcal{Z}^1(M)\}$  
has the property of finite intersections whenever $\psi$ is assumed to be sufficiently 
$C^0-$closed to the identity map. The conclusion follows from 
the Heine-Borel theorem from general topology. 
\end{proof}

 \begin{theorem}\label{Fix-Points3-5-7} $ $
 Let $(M,\Omega)$ be a closed oriented manifold. If $\psi\in \ker S_\Omega$ is  
sufficiently $C^0-$closed to the identity map, then
$$ \sharp\mathcal{FIX}(\psi) 
= \sharp\left(\bigcap_{\alpha\in \mathcal{Z}^1(M)}
\mathcal{N}^\varDelta(\psi,\alpha)\right)\geq 1.$$
\end{theorem}
\begin{proof} $ $
  Let $\psi\in \ker S_\Omega$ be 
sufficiently $C^0-$closed to the identity map. By Corollary \ref{Fix-Points3-5}, for each 
$p\in \bigcap_{\alpha\in \mathcal{Z}^1(M)}\mathcal{N}^\varDelta(\psi,\alpha)$, the 
pair $(\psi,p)$ is a zero for $\varDelta$, i.e., $\varDelta(\psi, \alpha)_p = 0$, 
for all $\alpha\in \mathcal{Z}^1(M)$. Thus, it follows from Theorem \ref{Functor3} 
that $p$ is a fix-point for $\psi$.
\end{proof}

The following result shows that fix-points can determine the flux 
geometry of certain volume-preserving 
diffeomorphisms.

\begin{lemma}\label{Rigiditylem2}{\bf (Weak $C^0-$flux conjecture I)} $ $
 Let $(M,\Omega)$ be a closed oriented manifold. Let $\{\psi_i\}\subset\ker{S}_\Omega$ be 
a sequence that converges in the $C^0-$metric to 
$\phi\in G_\Omega(M)$ such that there exists $\Phi\in Iso(\phi)_\Omega$ which satisfies the metric condition 
$d_{C^0}(\Phi, \Psi_i)\textless r(g)$ 
for some vanishing-flux path $\Psi_i\in  Iso(\psi_i)_\Omega$ whenever $i$ is sufficiently 
large. 
Then, we must have  $\phi\in \ker S_\Omega$. 
\end{lemma}
\begin{proof}$ $
 Let $\{\psi_i\}$ be a sequence of vanishing-flux volume-preserving 
diffeomorphisms which $C^0-$converges  
to a volume-preserving diffeomorphism $\phi$. Let $z_i\in \mathcal{FIX}(\psi_i)$, and derive 
that 
$\varDelta(\phi^i,\alpha)_{z_i}= 0 $, for each $i$, and for each closed $1-$form $\alpha$ 
 because $\mathcal{O}_{z_i}^{\Psi_i} $ is contractible (by Lemma \ref{FixPt4}). 
So, it follows 
from Lemma \ref{Cont1} that we must have 
\begin{equation}\label{Rigidityflux5}
  0= \lim_{d}\left( \lim_{d_{C^0}}\varDelta(\psi_i,\alpha)_{z_i}\right) =
 \lim_{d}\left(\varDelta(\phi,\alpha)_{z_i}\right),
\end{equation}
and hence, (\ref{Rigidityflux5}) together with Corollary \ref{GF10} 
imply that,  
\begin{equation}\label{Rigidityflux06}
\langle P(\alpha),\widetilde{S}_\Omega(\Phi) \rangle = 
Vol_\Omega(M)\lim_{d}\int_{\mathcal{O}_{z_i}^{\Phi}}\alpha,
\end{equation}
 for each closed $1-$form $\alpha$. On the other hand, since by assumption 
we have $d_{C^0}(\Phi, \Psi_i)\textless r(g)$ whenever $i$ is sufficiently 
large, we fix a positive integer $j$ large enough and derive that for each fixed $t$ and each $i\geq j$, 
one can connect $\Phi^t(z_i)$ to $\Psi^t_i(z_i)$ 
via a unique minimal geodesic $\xi_i^t$: $\mathcal{O}_{z_i}^{\Phi}$, $\mathcal{O}_{z_i}^{\Psi_i}$, and 
$ \xi_i^1$ delimit a $2-$chain. Hence, Stokes' theorem implies that for each $i\geq j$, we have  
$$\int_{\mathcal{O}_{z_i}^{\Phi}}\alpha = \int_{\mathcal{O}_{z_i}^{\Psi_i}}\alpha - 
\int_{\xi_i^1}\alpha = - \int_{\xi_i^1}\alpha,$$ for all closed $1-$form $\alpha$ 
because $\mathcal{O}_{z_i}^{\Psi_i}$ is contractible. That is, (\ref{Rigidityflux06}) implies that
$$
|\langle P(\alpha),\widetilde{S}_\Omega(\Phi) \rangle| = Vol_\Omega(M)
\lim_{d}|\int_{\xi_i^1}\alpha|\leq Vol_\Omega(M)
\|\alpha\|_{L^2}d(\phi,\psi_i)\rightarrow0, i\rightarrow\infty.
$$
i.e., 
$ \langle P(\alpha),\widetilde{S}_\Omega(\Phi) \rangle =  0$, 
for each closed $1-$form $\alpha$, i.e., $\widetilde{S}_\Omega(\Phi) = 0$: 
This implies that $\phi\in\ker_\Omega$.
\end{proof}

The following result generalizes and solves a weak version of the so-called $C^0-$flux conjecture.
 \begin{lemma}\label{Rigiditylem3}{\bf (Weak $C^0-$flux conjecture II)} $ $
 Let $(M,\Omega)$ be a closed oriented manifold. Let $\{\psi_i\}\subset\ker{S}_\Omega$ be 
a sequence that $C^0-$converges  to 
$\phi\in G_\Omega(M)$ such that:
\begin{itemize}
\item $\mathcal{ FIX}(\phi)\neq \emptyset$, and there exists 
$p \in \mathcal{ FIX}(\phi)$ such that $\mathcal{O}_{p}^{\Phi}$ 
is contractible for  some $\Phi\in Iso(\phi)_\Omega$, 
\item there exists a sequence $\{z_i\}\subset M$ with $z_i\in \mathcal{FIX}(\psi_i)$ 
that converges to $p$.
\end{itemize}
Then, we must have  $\phi\in \ker S_\Omega$. 
\end{lemma}
\begin{proof} $ $
  By assumption, since there exists a sequence 
$\{z_i\}\subset M$ with $z_i\in \mathcal{FIX}(\psi_i)$ 
that converges to $p$, then it follows 
from Lemma \ref{Cont1} that we must have\\ 
$
 0 = \|\alpha\|_{L^2}\lim_{d}\left(\lim_{d_{C^0}} \varDelta(\psi_i,\alpha)_{z_i}\right) 
=\|\alpha\|_{L^2}\varDelta(\phi,\alpha)_{p} = \langle P(\alpha),\widetilde{S}_\Omega(\Phi) \rangle - 0,
$
for each non-trivial closed $1-$form $\alpha$ because $\mathcal{O}_{p}^{\Phi}$ 
is contractible i.e. $
\langle P(\alpha),\widetilde{S}_\Omega(\Phi) \rangle = 0,$ 
 for each  closed $1-$form $\alpha$. Therefore, we have that $\widetilde{S}_\Omega(\Phi) = 0$. 
In particular, this 
implies that $\phi\in\ker_\Omega$
\end{proof}

\subsection{$C^0-$rigidity}
Let $\{\Phi_i =(\phi_i^t)\}$ be a sequence of vanishing-flux volume-preserving isotopies 
which converges uniformly to a path $\Phi=\{\phi_t\}$. Assume that for each $i$, the 
time-one map $\phi_i^1$ has a fix-point $z_i$.
It is known that $\phi$ has at least one fix-point: Assume that $d(x,\phi(x))\neq 0,$ for all 
$x\in M$, which implies that\\ 
$0\textless l_0:=\inf_{x\in M}d(x,\phi(x))\leq d(\phi(z_i),z_i)
\leq d(\phi(z_i),\phi(z_{i})\leq d_{C^0}(\phi_i,\phi) \rightarrow 0, 
i\rightarrow\infty.$
Thus, there exists $z\in M$ such that $\phi(z) = z$. 
In general, we have no information about the contractibility of the loop $\mathcal{O}_z^\Phi$
 (no matter the fact that all the orbits $ \mathcal{O}_{z_i}^{\Phi_i}$ are contractible): 
\begin{itemize}
 \item If $\Phi$ is smooth and belongs to the kernel of the flux, then 
$\mathcal{O}_z^\Phi$ is contractible (by Lemma \ref{FixPt4}).
\end{itemize}
In this subsection,  we will see that:
\begin{itemize}
 \item If the above $\phi$ is the identity map, then 
the orbit $\mathcal{O}_z^\Phi$ contractible is contractible whenever it is smooth.
\item  Furthermore, another result of the following section will show that if $\Phi$ is smooth, then 
$\mathcal{O}_z^\Phi$ is contractible.
\end{itemize}

Here are some rigidity results. 

\begin{theorem}\label{Rigidity} $ $
 Let $(M,\Omega)$ be a closed oriented manifold. If $\{\Phi_i\}$ is a sequence of 
vanishing-flux isotopies which $C^0-$converges to a loop $\gamma$, then 
any smooth orbit of $\gamma$  is contractible.
\end{theorem}
\begin{proof} $ $
 Let $z\in M$ be a point whose orbit under $\gamma$ is smooth, i.e. $\mathcal{O}_z^{\gamma}$ 
is smooth. By definition of $\gamma := (\gamma_t)_{t\in[0,1]}$, there exists $\{\Phi_i = (\phi_i^t)\}$ 
 a sequence of vanishing 
flux isotopies which converges uniformly to the loop $\gamma$. Since $\gamma_1 = id_M$, 
then the sequence of time-one 
maps $\{\phi_i^1\}$  
converges uniformly to the identity map $id_M$, and then it follows from  
Lemma \ref{Cont1} that $\lim_{C^0}\varDelta( \phi_i^1,\alpha)_z = \varDelta( id_M,\alpha)_z = 0$. 
On the other hand, by Corollary \ref{GF10}, we also have  
$ \|\alpha\|_{L^2}\varDelta( \phi_i^1,\alpha)_z =  -Vol_\Omega(M)
\int_{\mathcal{O}_{z}^{\Phi_i}}\alpha,$ 
for each $i$, and for all 
non-trivial closed $1-$form  $\alpha$. Hence, combining the above arguments give
 $
 \lim_{C^0}\left( \int_{\mathcal{O}_z^{\Phi_i}}\alpha\right) = 0,$
for all non-trivial closed $1-$form  $\alpha$. Observe that,
$
 |\int_{\mathcal{O}_z^{\gamma}}\alpha|
 \leq |\int_{\mathcal{O}_z^{\Phi_i}}\alpha -  \int_{\mathcal{O}_z^{\gamma}}\alpha| 
+ |\int_{\mathcal{O}_z^{\Phi_i}}\alpha|,$ 
for all $i$, and derive that $\int_{\mathcal{O}_z^{\gamma}}\alpha = 0,$ because 
$\lim_{C^0}\left( \int_{\mathcal{O}_z^{\Phi_i}}\alpha\right) = 0,$ 
and as in the proof of Lemma \ref{Cont1}, for $i$ sufficiently large, one 
can construct a smooth $2-$chain whose boundary consists of $\mathcal{O}_z^{\gamma}$, $ \mathcal{O}_z^{\Phi_i}$ 
and a 
minimizing geodesic $\chi_z^i$ that connects $ \gamma_1(z)$ to $\phi_i^1(z)$; and derive from Stokes' theorem that 
$
 |\int_{\mathcal{O}_z^{\Phi_i}}\alpha -  \int_{\mathcal{O}_z^{\gamma}}\alpha| 
= |\int_{\chi_z^i}\alpha|\leq \|\alpha\|_{L^2} d_{C^0}(\gamma, \Phi_i)\rightarrow0, i\rightarrow \infty.
$
Therefore, we have proved that $\int_{\mathcal{O}_z^{\gamma}}\alpha = 0,$ for all 
closed $1-$form $\alpha$, i.e., the $1-$cycle $ \mathcal{O}_z^{\gamma}$ is a boundary. 
\end{proof}

Theorem \ref{Rigidity} implies that a continuous Hamiltonian loop 
 has the property that: Any of its orbits which 
is smooth is contractible.\\

Here is a separation result.

 \begin{lemma}\label{Separation} $ $
 Let $(M,\Omega)$ be a closed oriented manifold. Let $\Phi$ be a volume-preserving isotopy out of 
$ \ker\widetilde{S}_\Omega$ with $\Phi(1) = \phi\neq id_M$. There exists $\delta_0\textgreater 0$
 which depends on $\Phi$ such that 
if the distance (in the $C^0-$topology) between $\Phi$ and the constant path identity is less than $\delta_0$, then 
none of the orbits of $\Phi$ is a minimal geodesic between its endpoints. 
\end{lemma}
\begin{proof} $ $
  Let $\Phi$ be a volume-preserving isotopy out of 
$ \ker\widetilde{S}_\Omega$, and let $\phi$ denote its time-one map. Since $\widetilde{S}_\Omega(\Phi)\neq0$, 
and $\langle,\rangle$ 
is non-degenerate, we can 
choose  a harmonic $1-$form $\mathcal{H}_0$ such that 
$\langle P(\mathcal{H}_0),\widetilde{S}_\Omega(\Phi) \rangle\neq 0$. Assume this done.
Let $r$ be the injectivity radius of the Riemannian metric on $M$, and suppose that $ \delta_0 = 
\dfrac{1}{\beta_0}\min\{r, \dfrac{|\langle P(\mathcal{H}_0),\widetilde{S}_\Omega(\Phi) \rangle|}
{\|\mathcal{H}_0\|_{L^2}Vol_\Omega(M)}\},$ for 
some $\beta_0\in]3,+\infty[$. 
If $d_{C^0}(\Phi, Id)\leq \delta_0,$ then one deduces as in the proof of 
Lemma \ref{Cont1} that for all $x\in M,$ we have 
\begin{equation}\label{SP1}
|\varDelta(\phi, \mathcal{H}_0)_x- 0| = |\varDelta(\phi, \mathcal{H}_0)_x
-\varDelta(Id, \mathcal{H}_0)_x| \leq 2Vol_\Omega(M)d_{C^0}(\Phi, Id),
\end{equation}
since $d_{C^0}(\phi, Id_M)\leq d_{C^0}(\Phi, Id)\textless r.$ 
On the other hand, derive from (\ref{SP1}) and Corollary \ref{GF10} 
that 
\begin{equation}\label{SP2}
\dfrac{|\langle P(\mathcal{H}_0),\widetilde{S}_\Omega(\Phi) \rangle|}
{Vol_\Omega(M)\|\mathcal{H}_0\|_{L^2}}\leq \dfrac{1}{Vol_\Omega(M)}|\varDelta(\phi, \mathcal{H}_0)_x| + 
\dfrac{1}{\|\mathcal{H}_0\|_{L^2}}|\int_{\mathcal{O}_{x}^\Phi}\mathcal{H}_0|,
\end{equation}
$$\leq 2d_{C^0}(\Phi, Id) + length(\mathcal{O}_{x}^\Phi),$$
for  each $x\in M$. Now, if $z\in M$ is such that $ length(\mathcal{O}_{z}^\Phi)$ 
minimizes the distance between its endpoints, then 
$ length(\mathcal{O}_{z}^\Phi)\leq d_{C^0}(\Phi, Id),$ and 
this together with (\ref{SP2}) implies that,\\ 
$
\beta_0|\langle P(\mathcal{H}_0),\widetilde{S}_\Omega(\Phi) \rangle|
\textless 3 |\langle P(\mathcal{H}_0),\widetilde{S}_\Omega(\Phi) \rangle|
$: This is a contradiction.
\end{proof}

\section{Metric geometries for volume-preserving diffeomorphisms}\label{Hofer-like} 
The goal of this section is to introduce and study right-invariant metrics on $G_\Omega(M)$. 
\subsection{A pseudo-norm on the group $G_\Omega(M)$}\label{Fixpoints-metric}
\subsubsection{A right-invariant pseudo-metric  on $G_\Omega(M)$}
For each $\psi\in G_\Omega(M)$, set 
\begin{equation}\label{Fixpoints-metric1}
 \rVert\psi\lVert^\infty := 
\sup_{\alpha\in \overline{\mathcal{B}^1_\mathcal{S}(0,1)}}
\left( \sup_{z\in M}\rvert\varDelta(\psi, \alpha)_z\lvert\right). 
\end{equation}

We have the following facts.

 \begin{proposition}\label{Fixpoints-metric2}
 The rule $ \rVert,\lVert^\infty$ has the following 
properties:
\begin{enumerate}
 \item Positivity: $\rVert\psi\lVert^\infty\geq 0$, for all $\psi\in G_\Omega(M)$.
\item Triangle inequality: 
$\rVert\psi\circ\phi\lVert^\infty\leq \rVert\psi\lVert^\infty + 
\rVert\phi\lVert^\infty,$ for all $\psi,\phi\in G_\Omega(M)$.
\item Symmetry: $\rVert\psi^{-1}\lVert^\infty = \rVert\psi\lVert^\infty$, for all $\psi\in G_\Omega(M)$. 
\item Pseudo-Nondegeneracy:  If $\rVert\psi\lVert^\infty= 0$, then  
$\langle P(\alpha),\widetilde{S}_\Omega(\Psi) \rangle =
Vol_\Omega(M)\int_{\mathcal{O}_{x}^\Psi}\alpha,$ for all $x\in M$,for all 
$\Psi\in Iso(\psi)_\Omega$, and
 for each $\alpha \in \mathcal{H}^1(M, \mathcal{S})$. 
\item Right-Invariant metric: The following rule 
\begin{equation}
 d^\Omega(\phi, \psi) := \rVert\phi\circ\psi^{-1}\lVert^\infty,
\end{equation}
for all $\psi,\phi\in G_\Omega(M)$, induces a right-invariant metric on $ G_\Omega(M)$.
\end{enumerate}

\end{proposition}
{\it Proof~:} 
\begin{enumerate}
 \item Positivity: This follows from the definition of  $ \rVert,\lVert^\infty$.
\item Triangle inequality: Let $\psi,\phi\in G_\Omega(M)$,
$\alpha\in \overline{\mathcal{B}^1_\mathcal{S}(0,1)}$, and derive with the help of Proposition \ref{GF1} that:
$$ \lvert\varDelta(\psi\circ\phi, \alpha)_z\rvert = 
\lvert  \varDelta(\phi, \alpha)_z + \varDelta(\psi, \alpha)_{\phi(z)}\rvert\leq 
\lvert  \varDelta(\phi, \alpha)_z\rvert + \lvert\varDelta(\psi, \alpha)_{\phi(z)}\rvert,
$$
for all $z\in M$, which implies the desired inequality.
\item Symmetry: From Proposition \ref{GF1}, we derive that 
$$0= \varDelta(id_M, \alpha) = \varDelta(\psi\circ\psi^{-1}, \alpha) 
= \varDelta(\psi^{-1}, \alpha) + \varDelta(\psi, \alpha)\circ\psi^{-1},$$ 
i.e.,  $ \varDelta(\psi^{-1}, \alpha) = - \varDelta(\psi, \alpha)\circ\psi^{-1}$. This implies that the 
desired equality.
\item Pseudo-Nondegeneracy:  If $\rVert\psi\lVert^\infty= 0$, then $ \varDelta(\psi, \alpha) = 0$, 
for all $ \alpha \in  \mathcal{H}^1(M, \mathcal{S})$. Thus, with the above vanishing condition, it follows 
from Corollary \ref{GF10} that\\ $
 \langle P(\alpha),\widetilde{S}_\Omega(\Psi) \rangle =
Vol_\Omega(M)\int_{\mathcal{O}_{x}^\Psi}\alpha,$ 
for all $x\in M$, for all $ \alpha \in  \mathcal{H}^1(M, \mathcal{S})$, and  
for each
$\Psi\in Iso(\psi)_\Omega$. 

\item  Right-Invariant metric:  We have,\\ 
 $d^\Omega(\phi\circ\rho, \psi\circ\rho) = \rVert(\phi\circ\rho)\circ(\psi\circ\rho)^{-1}\lVert^\infty
= \rVert \phi\circ\rho\circ\rho^{-1} \circ\psi^{-1}\lVert^\infty 
 = d^\Omega(\phi, \psi),$ 
for all $\rho, \psi,\phi\in G_\Omega(M)$.
\end{enumerate}
Q.E.D $\bigstar$\\

{\bf Conjecture:} The restriction $\rVert ,\lVert^\infty_{\ker}$ of the pseudo-norm 
$\rVert ,\lVert^\infty$ to the sub-group $\ker S_\Omega$ is non-degenerate.\\

Here is a consequence of the proof of Lemma \ref{Cont1}.
\begin{corollary}
 If $\psi\in \ker S_\Omega$ is sufficiently $C^0-$closed to the identity map,\\ then 
$\rVert \psi\lVert^\infty_{\ker}\leq 2Vol_\Omega(M)d_{C^0}(\psi, id_M)$. 
\end{corollary}
\begin{definition}
 For each non-empty open subset $\mathcal{U}\subset M$, 
we define the displacement energy of  $\mathcal{U}$ as:
$E^\Omega(\mathcal{U}):= \inf\{\rVert\psi\lVert^\infty: 
\psi\in G_\Omega(M), \psi(\mathcal{U}) \cap \mathcal{U} = \emptyset\}$. 
\end{definition}
We have not yet verified the positivity of $E^\Omega(\mathcal{U}) $  
whenever $\mathcal{U}$ is non-empty.\\

 {\bf Fact:} If it happens that $E^\Omega(\mathcal{U})\textgreater 0 $,  
whenever $\mathcal{U}$ is non-empty, then the following fact holds: Let 
$\{\phi_i\}\subset G_\Omega(M)$, let $\phi\in G_\Omega(M)$, and let $\psi: M\rightarrow M$ be any 
map such that $\{\phi_i\}\xrightarrow{C^0}\psi$, and $\rVert\phi_i\circ\phi^{-1}\lVert^\infty\rightarrow0, i
\rightarrow\infty $. Then, we must have $\psi =\phi$.\\

This could be used to define a new class of volume preserving-homeomorphisms. 
\subsection{Hofer-like geometry}
Assume that $(M, \omega)$ is a  closed symplectic manifold. The 
goal here is to  apply 
some results found in the previous sections in the study of 
 Hofer-like geometry of the group $Ham(M,\omega)$ of all Hamiltonian diffeomorphisms. 
As in \cite{TD2}, we will use the notation of Subsection \ref{De Rham}, with $\ast = 1$.
For each $\alpha\in \mathcal{Z}^1(M)$ we shall 
call the $1-$form $(\alpha - \mathcal{S}(P(\alpha)))$ the exact part of $\alpha$, and throughout 
all the paper, for simplicity, when this will be necessary, the latter $1-$form
 will be denoted $df_{\alpha,\mathcal{S} }$ 
to mean that it is the differential of a certain function 
that depends on $\alpha$ and $\mathcal{S}$; while we shall 
call the $1-$form $\mathcal{S}(P(\alpha))$ the $\mathcal{S}-$form of $\alpha$.
 Let $\mathbb H^1(M,\mathcal{S})$ denote the space of all $\mathcal{S}-$forms, and 
$\mathbb{B}^1(M)$ denote the set $(\mathcal{Z}^1(M) \backslash \mathbb H^1(M,\mathcal{S}))\cup\{0\}$.  
Consider $ \nu^B$ to be any norm on $\mathbb{B}^1(M)$. 
From 
now on, 
we equip $\mathbb H^1(M,\mathcal{S})$ with the $L^2-$norm $\lVert.\rVert_{L^2}-$Hodge norm for closed $1-$forms, 
and for  each 
non-negative $\lambda,$ we define a  semi-norm $N_{\mathcal{S}, \lambda}$ on  $\mathcal{Z}^1(M) $ as follows:
$
 N_{\mathcal{S}, \lambda}(\alpha) = \nu^B((\alpha - \mathcal{S}(P(\alpha)))) 
+ \lambda\lVert\mathcal{S}(P(\alpha))\rVert_{L^2},
$
for all $\alpha\in \mathcal{Z}^1(M) $.
\subsection{Hofer-like lengths}
In \cite{TD2}, given a symplectic isotopy $\Phi$ generated by $ (U,\mathcal{H}),$ 
the $L^{(1,\infty)}-$version and the $L^{\infty}-$version of Hofer-like lengths of a $\Phi$  
 are defined respectively  by
\begin{equation}\label{blg1}
 l_{\lambda, \mathcal{S}}^{(1,\infty)}(\Phi) = \int_0^1 \left( \nu^B(dU_t) + 
\lambda\lVert\mathcal{H}_t\rVert_{L^2}\right) dt, 
\end{equation}
\begin{equation}\label{blg2}
 l^\infty_{\lambda, \mathcal{S}}(\Phi) 
= \max_{t\in[0,1]}( \nu^B(dU_t) + \lambda\lVert\mathcal{H}_t\rVert_{L^2}).
\end{equation}
In the case that $H^1(M,\mathbb{R})$ vanishes, the above lengths are called 
Hofer's lengths whenever $\nu^B$ is the oscillation norm, i.e., 
$\nu^B(dU_t):= \max_{z\in M}(U_t(z)) -  \min_{z\in M}(U_t(z)) = osc(U_t)$, for each $t$ 
\cite{Hof-Zen94}.   
\begin{remark}\label{Lengthcon}
 We always have
\begin{equation}
 l_{\lambda, \mathcal{S}}(\Phi\ast_r\Psi) \leq   (\tau(\Phi) + 1 )l_{\lambda, \mathcal{S}}(\Psi) 
+ l_{\lambda, \mathcal{S}}(\Phi),
\end{equation}

where $\tau(\Phi)$ is positive constant which 
depends on $\Phi$; and 
\begin{equation}
  l_{\lambda, \mathcal{S}}^{(1,\infty)}(\Phi\ast_l\Psi) = l_{\lambda, \mathcal{S}}^{(1,\infty)}
(\Psi) +  l_{\lambda, \mathcal{S}}^{(1,\infty)}(\Phi).
\end{equation}
One can choose the smooth function $f$ in Subsection $2.1.2$ with bounded speed at any point, 
and  assume that

\begin{equation}\label{Lengthcon1}
 l_{\lambda, \mathcal{S}}^\infty (\Phi\ast_l\Psi) \leq \zeta \left(l_{\lambda, \mathcal{S}}^\infty (\Psi) + 
l_{\lambda, \mathcal{S}}^\infty (\Phi)\right),
\end{equation}
for  some $\zeta\in]1,6/5]$. 
\end{remark}
\subsection{Hofer-like norms}
The  $L^{(1,\infty)}-$energy and the $L^{\infty}-$energy of $\phi\in Symp_0(M,\omega)$ are respectively defined by:
\begin{equation}
 e_{\lambda, \mathcal{S}}(\phi) = \inf(l_{\lambda, \mathcal{S}}^{(1,\infty)}(\Phi)),
\end{equation}

and,
\begin{equation}
 e^{\infty}_{\lambda, \mathcal{S}}(\phi) = \inf(l_{\lambda, \mathcal{S}}^{\infty}(\Phi)),
\end{equation}
where each infimum is taken over the set of all symplectic isotopies 
$\Phi$ with time-one map equal to $\phi$. 
Therefore, the $L^{(1,\infty)}-$version and the $L^{\infty}-$version 
of the Hofer-like norms of $\phi$ are respectively defined by, 
\begin{equation}
 \|\phi\|_{\lambda, \mathcal{S}}^{(1,\infty)} 
= (e_{\lambda, \mathcal{S}}(\phi) + e_{\lambda, \mathcal{S}}(\phi^{-1}))/2,
\end{equation}
and
\begin{equation}
 \|\phi\|_{\lambda, \mathcal{S}}^\infty 
= (e^{\infty}_{\lambda, \mathcal{S}}(\phi) + e^{\infty}_{\lambda, \mathcal{S}}(\phi^{-1}))/2.
\end{equation}

The study of the Hofer-like geometry of $Ham(M,\omega)$ can be supported by the fact that it may exist at 
least one Hamiltonian diffeomorphism which can be connected to the identity 
by a symplectic isotopy which is not completely contained in $Ham(M,\omega)$ and 
which realizes the infimum of Hofer-like length: Given $\psi\in Ham(M,\omega)$, and a Hamiltonian isotopy 
$\Phi$ from the identity, pick any loop $\Psi\in \pi_1(Symp_0(M,\omega))$ whose 
flux is non-trivial, then the isotopy $\Xi:= \Psi\circ \Phi$ connects the identity map to $\phi$, but the 
isotopy $\Xi$ escapes  $Ham(M,\omega)$ since $ \widetilde{S}_\omega(\Xi)= 0 + \widetilde{S}_\omega(\Psi) \neq 0$.
Otherwise, the  
Hofer-like norm restricted to $Ham(M,\omega)$ 
is the Hofer metric.\\ Here are some examples of 
non-Hamiltonian isotopies with trivial flux whose endpoints belong to $Ham(M,\omega)$: Let consider 
a suitable smooth positive non-trivial function $\theta$ from $[0,1]$ to $[0,1]$ that vanishes near  
$0$ and $1$. Let $\alpha$  
be any closed $1-$form such that $P(\alpha)\neq 0$, and let $\Psi$ be the symplectic flow generated by  $\alpha$. 
The reparameterized path $ \Phi:s\mapsto \Psi(\theta(s))$ 
is a non-Hamiltonian loop, and $\widetilde{S}_\omega(\Phi) = 0$. So, if
 $\psi\in Ham(M,\omega)\backslash\{id_M\}$, and $\xi\in Iso(\psi)_\omega$ 
be any Hamiltonian isotopy, then the path $\xi\ast_r\Phi$ is a non-Hamiltonian isotopy that 
belongs to  $Iso(\psi)_\omega,$  whose 
flux is trivial. 
The above illustrations seem to suggest that the Hofer-like geometry 
of the group $Ham(M,\omega)$ is not trivial because the usual Hofer length of $\xi\ast_r\Phi$ does not make sense 
while its Hofer-like length does.\\ 

Although, it is proved in \cite{BusLec11} that the Hofer-like metric restricted to 
$Ham(M,\omega)$ is equivalent to the usual Hofer metric. This only tell 
us that the induced topologies on $Ham(M,\omega)$ from both Hofer and Hofer-like metrics
 are almost the same, and not more than that.\\  
In fact on a same space, two metrics can be equivalent while they induce 
different geometries (so many illustrations of this fact can be found in Euclidean spaces). 
The above arguments
 drive our attention to question:\\
   
{\it How does the locus of non-Hamiltonian symplectic isotopies with endpoints in $Ham(M,\omega)$ influences  
the geometry of $Ham(M,\omega)$?} In other words, how to detect the effect of flux geometry 
on the study of the metric geometry of $Ham(M,\omega)$?\\

To give some possible answers of the above question, we shall need the following alternate 
proof of a result found in \cite{Ban78, McDuff-SAl}.
This will be useful in the computation of the 
Hofer-like lengths of some isotopies arising from some deformations.  

\begin{theorem}(\cite{Ban78, McDuff-SAl})\label{Fgeo} $ $ Let $\Phi$ be a symplectic isotopy.
 If  $\widetilde{S}_\omega(\Phi) = 0$, then $\Phi$ is 
homotopic relatively to fixed endpoints to a Hamiltonian isotopy. 
\end{theorem}

\begin{proof} $ $

Assume that $\Phi = (\phi_t)$ is generated by $(U,\mathcal{H})$ (see \cite{BanTch1, TD2}). By assumption we have 
$\widetilde{S}_\omega(\phi_{(0,\mathcal{H})}) = [\int_0^1\mathcal{H}_tdt] = 0.$ Thus, it follows from 
Hodge's theory that $\int_0^1\mathcal{H}_tdt = 0$ since $\int_0^1\mathcal{H}_tdt$ is harmonic and $M$ 
is compact \cite{FWa}. 
For all $t$, as in \cite{McDuff-SAl}, let $s\mapsto\theta^t_s$ be the flow generated by the symplectic 
vector field $Y_t = -\int_0^tX_udu$ where $\iota(X_u)\omega = \mathcal{H}_u $ for all $u\in [0,1]$,  i.e., 
\begin{equation}
 Y_t(\theta_s^t(x)) = \dfrac{d}{ds}(\theta_s^t(x)),
\end{equation}
for all $x\in M$ with $s\in \mathbb{R},$ 
and $\theta^t_0 = id_M$. 
Since $Y_1 = 0 = Y_0$, we derive that $ \theta^0_s = id_M = \theta^1_s$ for all $s$.  
The map $ H : [0, 1]\times [0, 1]\rightarrow Symp_0(M,\omega),$ $
 (s,t)\mapsto \phi_{(0,\mathcal{H})}^t\circ\theta^t_s,$ induces a homotopy between 
the paths $t\mapsto \varphi_t = \phi_{(0,\mathcal{H})}^t\circ\theta^t_1$ 
and $ t\mapsto\phi_{(0,\mathcal{H})}^t$. On the other hand, if $V_{s,t}$ is the $2-$parameters 
family of symplectic vector fields 
defined as  
\begin{equation}
 V_{s,t}(\theta_s^t(x)) = \dfrac{\partial}{\partial t}(\theta_s^t(x))
\end{equation}
for all $x\in M$, 
then it follows from \cite{Ban78}-Proposition II.3.1.
 that 
\begin{equation}
 \dfrac{\partial}{\partial t}(Y_t) - \dfrac{\partial}{\partial s}(V_{s,t}) = [Y_t,V_{s,t}],
\end{equation}
for all $s, t$, i.e.,
\begin{equation}
 X_t - \dfrac{\partial}{\partial s}(V_{s,t}) = [Y_t,V_{s,t}].
\end{equation}
Integrating the latter equality in  the variable $s\in [0,u]$ leads to 
\begin{equation}
 \iota(V_{u,t})\omega= -u\iota(X_t)\omega  + df_{(u,t)},
\end{equation}
for all $(t,u)\in [0, 1]\times [0, 1]$, where  the function $x\mapsto f_{(u,t)}(x)$ 
is the normalized function  of the function  
$x\mapsto\int_0^t\int_0^u\omega(X_a, V_{s,t})(x)dsda,$ 
for all $(t,u)\in [0, 1]\times [0, 1]$ (see \cite{Ban10}). Hence, for 
each $u$, the isotopy $t\mapsto \theta_u^t$  is generated by $(f_{(u,.)}, -u\mathcal{H}).$ 
Thus, it follows that the generator of the isotopy 
$t\mapsto \varphi_t$  is given by
\begin{equation}
 (0, \mathcal{H})\Join (f_{(1,.)}, -\mathcal{H}) = (f_{(1,.)}\circ\phi_{(0,\mathcal{H})}^{-1} - 
\widetilde{\mathcal{F}_{\mathcal{H}_t}(\phi_{(0,\mathcal{H})}^{-1})(t)}, 0 ),
\end{equation}
where $\widetilde{\mathcal{F}_{\mathcal{H}_t}(\phi_{(0,\mathcal{H})}^{-1})(t)}$ is the normalized 
function of the function $\mathcal{F}_{\mathcal{H}_t}(\phi_{(0,\mathcal{H})}^{-1})(t)$ (see \cite{BanTch1, TD2}).  
This shows that the isotopy $t\mapsto \varphi_t$ is of the form $\phi_{(V,0)}$; hence  
the isotopy $t\mapsto \varphi_t\circ\phi_{(U,0)}^t$ is Hamiltonian, and is homotopic 
relatively to fixed endpoints to $\Phi$. 
\end{proof}
\begin{remark}
$ $ 
 In the rest of this section, we will assume that the norm $\nu^B$ is equivalent to the usual oscillation 
norm, and satisfies (if necessary) the property 
$
 \nu^B(\rho^\ast (df)) = \nu^B(df),
$
 for each 
smooth function $f$, and all $\rho\in Symp_0(M,\omega)$. 
\end{remark}
\begin{remark}\label{compari}  $ $ Given any symplectic vector field 
$X$, we can define its norm as follows:
\begin{equation}
 |X|_{\lambda, \mathcal{S}}:= \nu^B(( \iota(X)\omega - \mathcal{S}(P(\iota(X)\omega)))) 
+ \lambda\lVert\mathcal{S}(P(\iota(X)\omega))\rVert_{L^2}.
\end{equation}
Let  $h$ be a closed $1-$form such that $P(h)\neq 0$. We define the quantity  $\|\mathcal{S}(P(h))\|_0  $ as 
follows: Let 
$\chi(M,\omega)$ denote the space of all symplectic vector fields and consider the following set: 
$\mathcal{U}(M,\omega) := \{X\in \chi(M,\omega) :  |X|_{\lambda, \mathcal{S}}= 1\}.$ Therefore, put:
\begin{equation}
 \|\mathcal{S}(P(h))\|_0 
:= \sup_{z\in M}\left(\sup_{X\in \mathcal{U}(M,\omega) }(|(\mathcal{S}(P(h)))(X)(z)|)\right).
\end{equation}
Furthermore, since the space $\mathbb H^1(M,\mathcal{S})$ is a finite 
dimensional vector space whose dimension $b_1(M)$ is the first Betti number of $M$. Let   
 $\{h_i^\mathcal{S}\}$  be  any fixed  basis on the space of all $\mathcal{S}-$forms 
such that  
 $\mathcal{S}(P(h)) = \sum_i\lambda_ih_i^\mathcal{S}$. If we define a norm of $\mathcal{S}(P(h)) $ as 
$|\mathcal{S}(P(h))|_1 = \sum_i|\lambda_i|,$ 
then it follows that $
  \|\mathcal{S}(P(h))\|_0\leq E_\mathcal{S}|\mathcal{S}(P(h))|_1,$ 
where 
\begin{equation}
 E_\mathcal{S}:= \max_{1\leq i\leq b_1(M)}\|h_i^\mathcal{S}\|_0.
\end{equation}
 Note that
 the quantity $E_\mathcal{S}$ just defined above may take a large value in certain cases. The main 
difficulty in this section will be to go round the possible 
infinite nature of the constant $E_\mathcal{S}$. $\maltese$
\end{remark}
\begin{remark}\label{Rmk0}  $ $ Let $\{\alpha_t\}$ be a smooth family of closed $1-$forms, and define 
a smooth family of symplectic vector fields $ \{X_t\}$ by setting:
\begin{equation}\label{RmkI}
 \iota(X_t)\omega =\mathcal{S}(P(\alpha_t))=:\mathcal{H}_t,
\end{equation}
for each $t$. Now, put
\begin{equation}\label{Rmk2}
 Z_{s,t} := tX_{st} -2s\int_0^{t}X_udu,
\end{equation}
for all $(s,t)\in [0,1]\times[0,1],$ and  construct a $2-$parameters 
family of symplectic diffeomorphisms $G_{s,t}$ by integrating $Z_{s,t}$ in the variable $s$. Also, 
we may define another $2-$parameters 
family of symplectic vector fields as:
\begin{equation}\label{Rmk3}
 V_{s,t}(G_{(s,t)}(x)) = \dfrac{\partial}{\partial t}G_{(s,t)}(x),
\end{equation}
 for all $(s,t)\in [0,1]\times[0,1],$ and $x\in M$. $\maltese$
\end{remark}


\begin{lemma}\label{lem1}
 Let $\{X_t\}$ and $\{V_{s,t}\} $ be the smooth families of vector fields defined in (\ref{RmkI}) and (\ref{Rmk3}). 
Then 
\begin{equation}
 \dfrac{\sup_{s,t}|V_{s,t}|_{\lambda, \mathcal{S}}}{ \max\{2;6 E_\mathcal{S}\}
(1 + \sup_{s,t}|V_{s,t}|_{\lambda, \mathcal{S}})} \leq \sup_{t}|X_t|_{\lambda, \mathcal{S}}.
\end{equation}
\end{lemma}
\begin{corollary}\label{co1} $ $
 Let  $\{X_t\}$, $ \{Z_{s,t}\}$ and $\{V_{s,t}\} $ be the smooth families of vector fields defined 
in (\ref{RmkI}), (\ref{Rmk2}) and (\ref{Rmk3}). Then, 
\begin{equation}
 osc(\int_0^u\omega(Z_{s,t},V_{s,t})ds)\leq 4 E_\mathcal{S} \left(\sup_{s,t}|V_{s,t}|_{\lambda, \mathcal{S}}\right)
\left(\sup_{t}|X_t|_{\lambda, \mathcal{S}}\right),
\end{equation}
for each $u$ and each $t$ fixed.
\end{corollary}
\begin{proof} $ $
 
 For each $u$, and each $t$ (fixed) compute 
\begin{equation}\label{EQ2}
 osc (\int_0^u\omega(Z_{s,t},V_{s,t})ds)\leq 2\int_0^u\sup_{x\in M}|\omega(Z_{s,t},V_{s,t})(x)|ds.
\end{equation}
For each fixed $s$, for all $x\in M$, derive by the mean of Remark \ref{compari} that 
\begin{equation}\label{EQ3}
 |\omega(Z_{s,t},V_{s,t})(x)|\leq |(\iota(Z_{s,t})\omega)
(V_{s,t})(x_{s,t})|\leq\| \iota(Z_{s,t})\omega\|_0|V_{s,t}|_{\lambda, \mathcal{S}} 
\end{equation}
$$\leq (t|\mathcal{H}_{st}|_1 + 2s\int_0^t|\mathcal{H}_{q}|_1dq)|V_{s,t}|_{\lambda, \mathcal{S}} E_\mathcal{S}  ,$$
where $\iota(X_t)\omega =: \mathcal{H}_t$ for all $t$. 
Hence, it follows from (\ref{EQ2}) and (\ref{EQ3}) that 
$$
 osc(\int_0^u\omega(Z_{s,t},V_{s,t})ds)\leq 2E_\mathcal{S}\int_0^u(t|\mathcal{H}_{st}|_1 
+ 2s\int_0^t|\mathcal{H}_{u}|_1du)|V_{s,t}|_{\lambda, \mathcal{S}} ds,$$
$$\leq  2E_\mathcal{S}\sup_{s,t}|V_{s,t}|_{\lambda, \mathcal{S}} \int_0^u(t|\mathcal{H}_{st}|_1 + 
2s\int_0^t|\mathcal{H}_{u}|_1du)ds$$ $$ 
\leq  2E_\mathcal{S}\sup_{s,t}|V_{s,t}|_{\lambda, \mathcal{S}}(ut\sup_{t}|X_t|_{\lambda, \mathcal{S}} + 
u^2 \sup_{t}|X_t|_{\lambda, \mathcal{S}})$$
$$\leq 4E_\mathcal{S}\sup_{s,t}\left(\sup_{s,t}|V_{s,t}|_{\lambda, \mathcal{S}}\right)
\left(\sup_{t}|X_t|_{\lambda, \mathcal{S}}\right).$$
\end{proof} 
\begin{proof} of Lemma \ref{lem1}$ $

According to  Proposition II.3.1. found in \cite{Ban78} we have\\ $
 \dfrac{\partial}{\partial s}(V_{s,t}) - \dfrac{\partial}{\partial t}(Z_{s,t}) = [V_{s,t}, Z_{s,t}],$ 
integrating this with respect to $s\in[0,u]$ yields\\ $
 V_{u,t} - \int_0^u\dfrac{\partial}{\partial t}(Z_{s,t})ds = \int_0^u[V_{s,t}, Z_{s,t}]ds,
$
for all $u\in[0,1]$ and all $t$, i.e.,\\ $
 V_{u,t} = uX_{tu} -u^2X_t + \int_0^u[V_{s,t}, Z_{s,t}]ds,$
for all $u\in[0,1]$ and all $t$, because 
\begin{equation}\label{EQ7}
 \int_0^u\dfrac{\partial}{\partial t}(Z_{s,t})ds= \int_0^u \dfrac{\partial}{\partial t}(tX_{st})ds - u^2X_t
= \dfrac{\partial}{\partial t}\left(\int_0^{tu} X_{a}da\right) - u^2X_t
=  uX_{tu} -u^2X_t.
\end{equation}
Combining (\ref{EQ7}) together with the fact that 
$\iota([V_{s,t}, Z_{s,t}])\omega = d\left(\omega(V_{s,t}, Z_{s,t})\right),$ we 
derive from Corollary \ref{co1} that for all $u\in[0,1]$ and all $t$, we have
\begin{equation}\label{EQ8}
 |V_{u,t}|_{\lambda, \mathcal{S}}\leq (u + u^2)\sup_{t}|X_t|_{\lambda, \mathcal{S}}, + 
4E_\mathcal{S}\sup_{s,t}|V_{s,t}|_{\lambda, \mathcal{S}}\sup_{t}|X_t|_{\lambda, \mathcal{S}} ,
\end{equation}
$$\leq \max\{2,4E_\mathcal{S}\}(1 + \sup_{s,t}|V_{s,t}|_{\lambda, \mathcal{S}}) 
\sup_{t}|X_t|_{\lambda, \mathcal{S}}.$$

\end{proof}

Here is an iterative property in Hofer-like geometry. 
\begin{proposition}\label{Positivelength} $ $
Let $(M,\omega)$ be a closed symplectic manifold, and $\Omega_0$ be the corresponding symplectic 
volume form. If  $\Psi$ is any symplectic isotopy which is not a loop 
such that $\widetilde{S}_{\Omega_0}(\Psi)\neq 0$, then
\begin{enumerate}

\item $l_{\lambda, \mathcal{S}}^{(1,\infty)}(\Psi^l)\rightarrow\infty$, as $l\rightarrow \infty$, 
\item $ 0\textless \dfrac{\lvert\langle P(\mathcal{H}_0), 
\widetilde{S}_{\Omega_0}(\Psi)\rangle\lvert}{ Vol_{\Omega_0}(M)}
\leq E_\mathcal{S}\left(\lim_{l\rightarrow\infty}
\left[\dfrac{l_{\lambda, \mathcal{S}}^{(1,\infty)}(\Psi^l)}{l}\right]\right)$,\\
for some non-trivial $\mathcal{S}-$form $\mathcal{H}_0$.
\end{enumerate}

\end{proposition}
\begin{proof} $ $

Let $\Phi =(\phi_t)$ be a symplectic isotopy such that $\widetilde{S}_{\Omega_0}(\Psi)\neq 0$, and 
then choose a suitable $\mathcal{S}-$form 
$ \alpha_0$ such that $\langle P(\alpha_0),\widetilde{S}_{\Omega_0}(\Psi) \rangle\neq 0.$ Hence, 
derive from Proposition \ref{SC41} that; 
\begin{equation}
 |\langle  P(\alpha_0), \widetilde{S}_{\Omega_0}(\Phi)\rangle| 
= |\int_M\mathcal{F}_{\alpha_0}(\Phi)(1)\Omega_0|\leq Vol_{\Omega_0}(M)|\int_0^1\alpha_0(\dot\phi_t)\circ\phi_t(z)dt|,
\end{equation}
for some $z\in M$ that realizes the supremum of the function $x\mapsto \mathcal{F}_{\alpha_0}(\Phi)(1)(x)$. That is, 
$|\langle  P(\mathcal{H}_0), 
\widetilde{S}_{\Omega_0}(\Phi)\rangle| \leq 
E_\mathcal{S} Vol_{\Omega_0}(M)l_{\lambda, \mathcal{S}}^{(1,\infty)}(\Phi),$
where $ \mathcal{H}_0 = \alpha_0/|\alpha_0|_1.$ 
\end{proof}


\begin{remark}\label{lenght-diameter} $ $
Proposition \ref{Positivelength} is particularly 
interesting in the sense that it tells us that: 
If $\phi\in\left(Symp_0(M,\omega) \backslash Ham(M,\omega)\right)$, then for all sufficiently 
large integer $l$, an isotopy of the form $\Phi^l$ where $\Phi$ is any element of $Iso(\phi)_\omega$ 
can minimize the length 
functional $l_{\lambda, \mathcal{S}}^{(1,\infty)}$ on $Iso(\phi^l)_\omega$,
 if and only if, the diameter of $Symp_0(M,\omega)$ 
with respect to the Hofer-like metric 
is sufficiently large (and even infinite).$ \maltese$
\end{remark}
\subsection{ Hofer-like geometry via loops}
The following result shows (with a certain restriction) the invariance 
of the Hofer-like energy with respect to concatenation by loops in $Symp_0(M,\omega)$, and 
 once more detects another effect of the flux geometry in the study of Hofer-like geometry.  
To see this, consider $\mathcal{L}oop(M,\omega)$ to be the set consisting all the loops 
at the identity in 
$Symp_0(M,\omega),$
pick $\phi\in Ham(M,\omega),$ and for each $\Phi \in Iso(\phi)_\omega,$  put 
\begin{equation}
 \Upsilon(\Phi) := \{\Psi\in \mathcal{L}oop(M,\omega): \widetilde{S}_\omega(\Phi) = \widetilde{S}_\omega(\Psi)\}.
\end{equation}
Since $\phi$ is Hamiltonian, then for each $\Phi \in Iso(\phi)_\omega,$ the set $\Upsilon(\Phi)$ is non-empty. This
 is supported by two facts found in \cite{Ban78}: The fist implies that $Ham(M,\omega)$ is path 
connected through smooth Hamiltonian isotopies, and the second implies that
 any path from the identity to $\phi$ has its flux in $\Gamma_\omega$. 

\begin{proposition}\label{proEQU}({\bf Energy invariance}). $ $
For each $\phi\in Ham(M,\omega),$ we have
$$ e_{\lambda, \mathcal{S}}(\phi) = 
\inf_{\Phi\in Iso(\phi)_\omega}\left
(\inf_{\Psi\in \Upsilon(\Phi)}(l_{\lambda, \mathcal{S}}^{(1,\infty)}(\Phi\ast_l\Psi))\right).$$
\end{proposition}
\begin{proof} $ $
 
 Pick $\phi\in Ham(M,\omega),$ and derive from remark (\ref{Lengthcon}) that 
\begin{equation}\label{proEQUcon}
 l_{\lambda, \mathcal{S}}^{(1,\infty)}(\Phi\ast_l\Psi) = l_{\lambda, \mathcal{S}}^{(1,\infty)}(\Psi) + 
l_{\lambda, \mathcal{S}}^{(1,\infty)}(\Phi),
\end{equation}
for all $\Phi\in Iso(\phi)_\omega$, and all $\Psi\in \Upsilon(\Phi)$. Then, in (\ref{proEQUcon}), 
for each $\Phi\in Iso(\phi)_\omega$ one 
passes to the infimum  over the set $ \Upsilon(\Phi)$ in a first time, 
and next passing to the infimum over $Iso(\phi)_\omega$ in a second time yields,
\begin{equation}\label{EQU-38}
 \inf_{\Phi\in Iso(\phi)_\omega}
\left(\inf_{\Psi\in \Upsilon(\Phi)}(l_{\lambda, \mathcal{S}}^{(1,\infty)}(\Phi\ast_l\Psi))\right) 
= \inf_{\Phi\in Iso(\phi)_\omega}(l_{\lambda, \mathcal{S}}^{(1,\infty)}(\Phi)) + 
\inf_{\Phi\in Iso(\phi)_\omega}\left(\inf_{\Psi\in \Upsilon(\Phi)}(l_{\lambda, \mathcal{S}}^{(1,\infty)}(\Psi))\right).
\end{equation}
But, $\phi$ being Hamiltonian, the set $Iso(\phi)_\omega$ contains  at least 
a Hamiltonian isotopy (see \cite{Ban78}), and so, when the isotopy $\Phi\in Iso(\phi)_\omega$ is Hamiltonian, 
then the trivial loop identity belongs 
to $\Upsilon(\Phi),$ i.e, at a certain time the loop $\Psi$ will take the value of the trivial loop identity 
so that  
$\inf_{\Phi\in Iso(\phi)_\omega}
\left(\inf_{\Psi\in \Upsilon(\Phi)}(l_{\lambda, \mathcal{S}}^{(1,\infty)}(\Psi))\right) = 0.$  
\end{proof}

A result found in \cite{TD2} states that the above energies agree everywhere. So, the above energy 
invariance result will continue to hold if we replace the energy functional   
$e_{\lambda, \mathcal{S}}$ by the energy functional $e_{\lambda, \mathcal{S}}^\infty$.\\


Here is an alternative proof of the non-degeneracy of the norm $\|.\|_{\lambda, \mathcal{S}}$. Another 
alternative proof of this result was given in \cite{TD2}. The difference 
between these alternate proofs come from the fact that in \cite{TD2} the given 
proof follows as a direct application of Theorem $3.3-$\cite{TD2}, but here  
we do not appeal to Theorem $3.3-$\cite{TD2}.

 \begin{lemma}\label{lem2} $ $
 If $\phi\in Symp_0(M,\omega)$ such that $\|\phi\|_{\lambda, \mathcal{S}}^\infty = 0$, then $\phi= id_M$.
\end{lemma}
\begin{proof} $ $

Let  $\phi\in Symp_0(M,\omega)$ such that $\|\phi\|_{\lambda, \mathcal{S}}^\infty = 0$, and as in \cite{Ban10}, we 
derive that 
$\phi\in Ham(M,\omega)$. The rest of the proof consists to show that the 
Hofer norm of $\phi$ is trivial. 
Since $e_{\lambda, \mathcal{S}}^\infty(\phi) = 0,$ we derive 
 that for all positive integer $N$ there is $\Phi_N\in Iso(\phi)_\omega$ such that:
\begin{equation}
 l_{\lambda, \mathcal{S}}^\infty(\Phi_N)\leq \exp(-N\max\{2; 4E_\mathcal{S}\}).
\end{equation}
We have $\widetilde{S}_\omega(\Phi_N) = 0$ (if necessary choose the $N$'s sufficiently large 
so that\\ $\exp(-N\max\{2; 4E_\mathcal{S}\})$ is arbitrarily small), because we always have 
\begin{equation}
 \|\widetilde{S}_\omega(\Phi_N) \|_{L^2} \leq l_{\lambda, \mathcal{S}}^\infty(\Phi_N),
\end{equation}
 and 
$\Phi_N$ having its endpoints in $Ham(M,\omega)$, must have its flux in the Calabi group which is discrete: 
So, when $N$ is sufficiently large, then $\|\widetilde{S}_\omega(\Phi_N) \|_{L^2}$ is arbitrarily small, 
and hence trivial.  
 Let $\rho_t^N\circ_\mathcal{S}\psi_t^N$ be the 
$\mathcal{S}-$decomposition of the isotopy $\Phi_N$ (see \cite{TD2}). Since 
$\widetilde{S}_\omega(\{\rho_t^N\}) = \widetilde{S}_\omega(\Phi_N) = 0,$ and
\begin{equation}
 \sup_{t}|V_{1,t}^N|_{\lambda, \mathcal{S}}\leq\sup_{t,s}|V_{s,t}^N|_{\lambda, \mathcal{S}},
\end{equation}
where the smooth $2-$parameter family of symplectic vector fields 
$\{V_{s,t}^N\}$ is constructed as in (\ref{Rmk3}) using 
$X_t^N := \dot\rho_t^N$ for each $t$, 
we derive from Lemma \ref{lem1} that 
\begin{equation}
 \dfrac{\|\rho_1^N\|_{Hofer} }{1 +\|\rho_1^N\|_{Hofer} }\leq \dfrac{\sup_{t}|V_{1,t}^N|_{\lambda, \mathcal{S}}}{(1 +
 \sup_{t}|V_{1,t}^N|_{\lambda, \mathcal{S}})}\leq \max\{2; 4E_\mathcal{S}\}\sup_t|X_t^N|_ {\lambda, \mathcal{S}}
\leq \max\{2; 4E_\mathcal{S}\}l_{\lambda, \mathcal{S}}^\infty(\Phi_N)
\end{equation}
 $$ \leq \max\{2; 4E_\mathcal{S}\}
\exp(-N\max\{2; 4E_\mathcal{S}\}).$$
 That is,
\begin{equation}\label{0infimum}
\dfrac{\|\rho_1^N\|_{Hofer}}{1 +\|\rho_1^N\|_{Hofer} }\leq\max\{2; 4E_\mathcal{S}\}\exp(-N\max\{2; 4E_\mathcal{S}\}).
\end{equation}
On the other hand, we also have 
\begin{equation}\label{infimum}
\|\psi_1^N\|_{Hofer} \leq l_{\lambda, \mathcal{S}}^\infty(\psi_t^N)\leq l_{\lambda, \mathcal{S}}^\infty
(\Phi_N) \leq \exp(-N\max\{2; 4E_\mathcal{S}\}).
\end{equation}
Thus, as $N\rightarrow\infty$, the right-hand side in (\ref{infimum}) (resp. the right-hand side in (\ref{0infimum})) 
tends to zero: Thus, we must have $\|\rho_1^N\|_{Hofer} \rightarrow 0, N\rightarrow\infty,$  
and $\|\psi_1^N\|_{Hofer} \rightarrow 0, N\rightarrow\infty$. Finally, since 
$\phi = \rho_1^N\circ\psi_1^N$, for each integer $N$, we derive from the triangle 
inequality that 
\begin{equation}
 \|\phi\|_{Hofer} = \|\rho_1^N\circ\psi_1^N\|_{Hofer}
\leq \|\rho_1^N\|_{Hofer} + \|\psi_1^N|_H\rightarrow 0, N\rightarrow\infty,
\end{equation}
 i.e., $ \|\phi\|_{Hofer} = 0$.  
Hence, from the non-degeneracy of the Hofer norm, we have $\phi = id_M$. 
\end{proof}

The main result of this section is the following comparison theorem. 

\begin{theorem}\label{EQU}  $ $
  Let $(M,\omega)$ be a closed symplectic manifold. Then, 
for all Hamiltonian diffeomorphism $\phi$ and for all positive $\epsilon$, 
there exist two elements $\rho^\epsilon$ and $\psi^\epsilon$ 
of $Ham(M,\omega)$  with $\phi = \rho^\epsilon\circ\psi^\epsilon$ 
such that $$
\|\phi\|_{Hofer} + \|\rho^\epsilon\|_{Hofer}\|\psi^\epsilon\|_{Hofer}\leq (1 
+\|\rho^\epsilon\|_{Hofer}) \left(\mathcal{G}(E_\mathcal{S})\|\phi\|_{\lambda, \mathcal{S}}^\infty + 
\epsilon\right),
$$
for some positive constant $\mathcal{G}(E_\mathcal{S})$ that depends on $E_\mathcal{S}$.
\end{theorem}

Here is a consequence of Theorem \ref{EQU}, this  agrees with a result found in \cite{BusLec11}.
\begin{proposition}\label{cP11} $ $
Let $(M,\omega)$ be a closed symplectic.
Therefore, there exists a positive constant $\kappa$ such that for all $\phi\in Ham(M,\omega)$, one has  
\begin{equation}
 \kappa\|\phi\|_{Hofer} \leq \|\phi\|_{\lambda, \mathcal{S}}^\infty.
\end{equation} 

\end{proposition}
\begin{proof} $ $

Let $ \{\phi_k\}\subset Ham(M,\omega)$ be a 
sequence that converges to the constant map identity 
in the Hofer-like norm. By
Theorem \ref{EQU}, for each positive integer $k \geq \nu$ 
(for some sufficiently large positive integer $\nu$, fixed) , and for $\epsilon = 1/k$, 
there exist two elements $\rho^\epsilon_k$ and 
$\psi^\epsilon_k$  of $Ham(M,\omega)$  with $\phi_k = \rho^\epsilon_k\circ\psi^\epsilon_k$ 
such that 
 \begin{equation}\label{constant1}
\dfrac{\|\phi_k\|_{Hofer} + \|\rho^\epsilon_k\|_{Hofer}\|\psi^\epsilon_k\|_{Hofer}}{1 
+\|\rho^\epsilon_k\|_{Hofer} }\leq \mathcal{G}(E_\mathcal{S})\|\phi_k\|_{HL} + 1/k.
\end{equation}
Since by assumption we have $ \|\phi_k\|_{\lambda, \mathcal{S}}^\infty\rightarrow 0, k\rightarrow\infty,$
then 
 assume that 
\begin{equation}
 \|\phi_k\|_{\lambda, \mathcal{S}}^\infty\leq \exp(-k \mathcal{G}(E_\mathcal{S})),
\end{equation}

for each  $k \geq \nu$ (if necessary modify $\nu$). That is, (\ref{constant1}) implies that 
\begin{equation}\label{constant3}
\dfrac{ \|\rho^\epsilon_k\|_{Hofer}\|\psi^\epsilon_k\|_{Hofer}}{1 +\|\rho^\epsilon_k\|_{Hofer} } + 
\dfrac{\|\phi_k\|_{Hofer} }{1 +\|\rho^\epsilon_k\|_{Hofer} }
\leq  \mathcal{G}(E_\mathcal{S})\exp(-k \mathcal{G}(E_\mathcal{S}))
 + 1/k,
\end{equation}
for each positive integer $k\geq \nu$. In particular, (\ref{constant3}) implies  that
\begin{equation}
 \dfrac{ \|\rho^\epsilon_k\|_{Hofer}\|\psi^\epsilon_k\|_{Hofer}}{1 
+\|\rho^\epsilon_k\|_{Hofer} } \rightarrow0, k\rightarrow\infty,
\end{equation}

 and 
\begin{equation}
  \dfrac{\|\phi_k\|_{Hofer} }{1 +\|\rho^\epsilon_k\|_{Hofer} }
\rightarrow 0,k\rightarrow\infty.
\end{equation}
So, the convergence $(\dfrac{\|\rho^\epsilon_k\|_{Hofer}}{1 
+\|\rho^\epsilon_k\|_{Hofer} })\|\psi^\epsilon_k\|_{Hofer}\rightarrow0, k\rightarrow\infty,$ implies that\\ either 
$(\dfrac{\|\rho^\epsilon_k\|_{Hofer}}{1 +\|\rho^\epsilon_k\|_{Hofer} })\rightarrow0, k\rightarrow\infty$, 
or $\|\psi^\epsilon_k\|_{Hofer}\rightarrow0, k\rightarrow\infty$.\\
{\bf Case (1):}  If 
\begin{equation}
 (\dfrac{\|\rho^\epsilon_k\|_{Hofer}}{1 +\|\rho^\epsilon_k\|_{Hofer} })\rightarrow0, k\rightarrow\infty,
\end{equation}
then  we derive that $\|\rho^\epsilon_k\|_{Hofer}\rightarrow0, k\rightarrow\infty;$ and this together 
with the convergence $\dfrac{\|\phi_k\|_{Hofer}}{1 +\|\rho^\epsilon_k\|_{Hofer} }\rightarrow 0,k\rightarrow\infty,$ 
implies that $\|\phi_k\|_{Hofer}\rightarrow0, k\rightarrow\infty$.\\ 
 {\bf Case (2):} If $\|\psi^\epsilon_k\|_{Hofer}\rightarrow0, k\rightarrow\infty,$ 
then from the fact the triangle inequality,
\begin{equation}
 \|\rho^\epsilon_k\|_{Hofer}\leq \|\psi^\epsilon_k\|_{Hofer} + \|\phi_k\|_{Hofer},
\end{equation} 
we derive that 
\begin{equation}
 \left(\dfrac{\|\phi_k\|_{Hofer}}{1 +\|\phi_k\|_{Hofer} + 
\|\psi^\epsilon_k\|_{Hofer} }\right)\leq\dfrac{\|\phi_k\|_{Hofer}}{1 +
\|\rho^\epsilon_k\|_{Hofer} }\rightarrow 0,k\rightarrow\infty.
\end{equation}
On the other hand, it follows from the continuity of the function $f_1: x\mapsto x/(1+x)$ that 
\begin{equation}\label{EQ13}
 f_1(\lim_{k\rightarrow\infty}(\|\psi^\epsilon_k\|_{Hofer} + \|\phi_k\|_{Hofer}))
= \lim_{k\rightarrow\infty}f_1(\|\psi^\epsilon_k\|_{Hofer} + \|\phi_k\|_{Hofer})   
\end{equation}
$$ = \lim_{k\rightarrow\infty}\left(\dfrac{\|\phi_k\|_{Hofer}}{1 +\|\phi_k\|_{Hofer} 
+ \|\psi^\epsilon_k\|_{Hofer} }\right) = 0.$$
Equation (\ref{EQ13}) implies that $\|\phi_k\|_{Hofer}\rightarrow0, k\rightarrow\infty$. 
Finally, we have proved that the norm $\|.\|_{\lambda, \mathcal{S}}^\infty$ 
restricted to $Ham(M,\omega)$, is equivalent to 
the Hofer norm $ \|.\|_{H}$. That is,
 there exists a positive constant $\kappa$ such that for all Hamiltonian diffeomorphisms 
$\phi$, one has $
 \kappa\|\phi|_{H}\leq \|\phi\|_{\lambda, \mathcal{S}}^\infty.
$ 
\end{proof}

\subsection{The ingredients of the proof of Theorem \ref{EQU}}\label{Proof-M} 
Recall that one can classify closed symplectic manifolds in two categories: one category 
consists of closed symplectic manifolds 
with trivial Calabi's groups, i.e., $\Gamma_\omega = \{0\}$, and the other consists of closed symplectic manifolds 
with  $\Gamma_\omega \neq \{0\}$.
Compact symplectic manifolds with vanishing flux groups includes oriented closed surfaces of genus bigger than
$1$. Furthermore, Kedra-Kotschick-Morita \cite{KKM} found a long list of compact
symplectic manifolds with vanishing flux groups.\\

We have the following fact. 
\begin{lemma}\label{Vanishingamma} $ $ Let $(M,\omega)$ be a closed symplectic manifold whose 
flux group $\Gamma_\omega$ is trivial. 
Then, for all Hamiltonian diffeomorphism $\phi$ and  for all positive $\epsilon$, 
there exist two elements $\rho^\epsilon$ and $\psi^\epsilon$ 
of $Ham(M,\omega)$  with $\phi = \rho^\epsilon\circ\psi^\epsilon$ 
such that 
\begin{equation}
 \|\phi\|_{Hofer} + \|\rho^\epsilon\|_{Hofer}\|\psi^\epsilon\|_{Hofer} \leq (1 +
\|\rho^\epsilon\|_{Hofer} ) \left(\mathcal{G}_1(E_\mathcal{S}) e_{\lambda, \mathcal{S}}^\infty(\phi) + \epsilon\right),
\end{equation}
for some positive constant $\mathcal{G}_1(E_\mathcal{S})$ which depends on $E_\mathcal{S} $. 
\end{lemma}
\begin{proof} $ $

Assume that the flux group $\Gamma_\omega$ is trivial. 
Let $\phi\in Ham(M,\omega),$ $\epsilon\textgreater0$, and consider $\Phi = \{\phi_t\}$ 
to be any symplectic isotopy from the identity to 
$\phi$.  Since  $\phi\in Ham(M,\omega),$ it follows from \cite{Ban78} 
that $\widetilde{S}_\omega(\Phi)\in\Gamma_\omega =\{0\}$, 
and hence 
we see that $\widetilde{S}_\omega(\Phi) = 0$. Now, set 
\begin{equation}
 \mathcal{G}_1(E_\mathcal{S}):= \max\{2; 4E_\mathcal{S}\},
\end{equation}
and derive from the characteristic of infimum that for $\tau = \epsilon/\mathcal{G}_1(E_\mathcal{S})$, 
 there exists an isotopy 
$\Phi_\tau\in Iso(\phi)_\omega$ such that $
 l_{\lambda, \mathcal{S}}^\infty(\Phi_\tau)\leq e_{\lambda, \mathcal{S}}^\infty(\phi) + \tau.
$
 Remember that 
if $\Phi_\tau$  has $\mathcal{S}-$decomposition $\rho_t^\tau\circ\psi_t^\tau$, then it follows from 
the lines of the proof of Lemma \ref{lem1} that,
$$
 \dfrac{\|\phi\|_{Hofer} + \|\rho_1^\tau\|_{Hofer}\|\psi^\tau_1\|_{Hofer}}{1 +
\|\rho^\tau_1\|_{Hofer} } \leq \dfrac{\|\rho^\tau_1\|_{Hofer}}{1 +\|\rho^\tau_1\|_{Hofer} } + 
l_{\lambda, \mathcal{S}}^\infty(\psi^\tau_t)$$ 
$$\leq
 \mathcal{G}_1(E_\mathcal{S})l_{\lambda, \mathcal{S}}^\infty(\rho_t^\tau) 
+ l_{\lambda, \mathcal{S}}^\infty(\psi^\tau_t)\leq 
\mathcal{G}_1(E_\mathcal{S})\left(l_{\lambda, \mathcal{S}}^\infty(\rho_t^\tau) 
+ l_{\lambda, \mathcal{S}}^\infty(\psi^\tau_t)\right),
$$
$$ =  \mathcal{G}_1(E_\mathcal{S})l_{\lambda, \mathcal{S}}^\infty(\Phi_\tau),$$
since $$\|\psi^\tau_1\|_{Hofer}\leq l_{\lambda, \mathcal{S}}^\infty(\psi^\tau_t),$$ 
$$\|\phi\|_{Hofer} = \|\rho_1^\tau\circ\psi^\tau_1\|_{Hofer} \leq  \|\rho_1^\tau\|_{Hofer} +  
\|\psi^\tau_1\|_{Hofer},$$ 
and $\mathcal{G}_1(E_\mathcal{S})\textgreater 1 $. That is,
\begin{equation}\label{infimum2}
\dfrac{\|\phi\|_{Hofer} + \|\rho_1^\tau\|_{Hofer}\|\psi^\tau_1\|_{Hofer}}{1 +\|\rho^\tau_1\|_{Hofer} }\leq 
\mathcal{G}_1(E_\mathcal{S})l_{\lambda, \mathcal{S}}^\infty(\Phi_\tau)\leq \mathcal{G}_1(E_\mathcal{S})
e_{\lambda, \mathcal{S}}^\infty(\phi) + \epsilon.
\end{equation}
\end{proof}

When the flux group $\Gamma_\omega$ is non-trivial, then for some 
$\phi\in Ham(M,\omega)$, the set $Iso(\phi)_\omega$ 
may contain a symplectic isotopy $\Xi$ whose flux is non-trivial: In this situation, the first item in   
Proposition \ref{Positivelength}  
is telling to us that, if $Iso(\phi)_\omega$ contains such an isotopy $\Xi$, 
 then iterating $\phi$ 
a large number of times, say, in $l-$times ($l$ sufficiently large) 
will yield $\phi^l\in Ham(M,\omega)$ while $\Xi^l\in Iso_\omega(\phi^l)$ 
may have  an infinite Hofer-like length. 
But, since the isotopy $\Xi$ has its flux in $\Gamma_\omega$, then 
one can always find $\Psi\in \mathcal{L}oop(M,\omega)$ 
such that the isotopy $\Psi^{-1}\ast_l\Xi$ has a trivial flux, whereas, 
both isotopies $\Psi^{-1}\ast_l\Xi$ and $\Xi$ 
still have the same endpoints.\\
In the sequel, we see that for more 
convenience in the control of the Hofer-like lengths, it could be judicious to 
kill the flux of isotopies by deforming them without changing their endpoints. However, such a 
deformation is only possible 
for any isotopy whose endpoints belong to $Ham(M,\omega)$. This suggests to 
us that the Calabi group is in control of the Hofer-like 
geometry of $Ham(M,\omega)$.This motivated the following result.

\begin{lemma}\label{Nonvanishingamma} $ $
 Let $(M,\omega)$ be a closed symplectic manifold whose flux group $\Gamma_\omega$ is non-trivial. 
Then, for all Hamiltonian diffeomorphism $\phi$, 
for all positive $\epsilon$, there exist two elements $\rho^\epsilon$ and $\theta^\epsilon$ 
of $Ham(M,\omega)$  with $\phi = \rho^\epsilon\circ\theta^\epsilon$ 
such that 
\begin{equation}
 \|\phi\|_{Hofer} + \|\rho^\epsilon\|_{Hofer}\|\theta^\epsilon\|_{Hofer} \leq (1 +
\|\rho^\epsilon\|_{Hofer})\left( \mathcal{G}_2(E_\mathcal{S})e_{\lambda, \mathcal{S}}^\infty(\phi) + \epsilon\right),
\end{equation}
for some positive constant $\mathcal{G}_2(E_\mathcal{S})$ which depends on $E_\mathcal{S} $.
\end{lemma}

\begin{proof} $ $

Assume that the flux group $\Gamma_\omega$ is non-trivial. Let $\phi\in Ham(M,\omega)$ 
and $\epsilon\textgreater0$. Since  
for each 
$\Phi\in Iso(\phi)_\omega$, we have that all  
 $\Psi\in \Upsilon(\Phi)$ satisfies $\widetilde{S}_\omega(\Phi) = \widetilde{S}_\omega(\Psi),$ then 
when $\Psi$ describes the set $\Upsilon(\Phi)$,  the fluxes of the symplectic 
paths $\Psi^{-1}\ast_l\Phi$ are trivial, and $\Psi^{-1}\ast_l\Phi\in Iso(\phi)_\omega.$ Since by 
Proposition \ref{proEQU} we have\\
$
 e_{\lambda, \mathcal{S}}^\infty(\phi)= 
\inf_{\Phi\in Iso(\phi)_\omega}
\left(\inf_{\Psi\in \Upsilon(\Phi)}(l_{\lambda, \mathcal{S}}^\infty(\Psi^{-1}\ast_l\Phi))\right), 
$
 then for $\tau = \epsilon/\mathcal{G}_1(E_\mathcal{S})$ 
where $\mathcal{G}_1(E_\mathcal{S})$ is the constant which appears in Lemma \ref{Vanishingamma}, 
there exist two isotopies 
$\Phi_\tau\in Iso(\phi)_\omega$  and $\Psi_\tau\in \Upsilon(\Phi_\tau)$ such that 
\begin{equation}\label{EQ16}
 l_{\lambda, \mathcal{S}}^\infty(\Psi_\tau^{-1}\ast_l\Phi_\tau)\leq \zeta 
e_{\lambda, \mathcal{S}}^\infty(\phi) + \tau,
\end{equation}
where $\zeta$ is the constant which appears in (\ref{Lengthcon1}). 
On the other hand, since $\widetilde{S}_\omega(\Psi^{-1}_\tau\ast_l\Phi_\tau) = 0 $, 
we derive as in the proof of Lemma \ref{Vanishingamma}  that  
if $\Psi^{-1}_\tau\ast_l\Phi_\tau$  has $\mathcal{S}-$decomposition $\rho_t^\tau\circ\theta_t^\tau$, then 
\begin{equation}\label{EQ17}
 \dfrac{\|\phi\|_{Hofer} + \|\rho_1^\tau\|_{Hofer}\|\theta^\tau_1\|_{Hofer}}{1 +\|\rho^\tau_1\|_{Hofer} }
\leq \mathcal{G}_1(E_\mathcal{S})l_{\lambda, \mathcal{S}}^\infty(\Psi^{-1}_\tau\ast_l\Phi_\tau),
\end{equation}
where $\mathcal{G}_1(E_\mathcal{S}) $ is the constant in Lemma \ref{Vanishingamma},  
i.e., combining (\ref{EQ16}) and (\ref{EQ17}) together yields, 
 \begin{equation}\label{infimum4}
\dfrac{\|\phi\|_{Hofer} + \|\rho_1^\tau\|_{Hofer}\|\theta^\tau_1\|_{Hofer}}{1 +\|\rho^\tau_1\|_{Hofer} }
 \leq \zeta\mathcal{G}_1(E_\mathcal{S})e_{\lambda, \mathcal{S}}^\infty(\phi) + \mathcal{G}_1(E_\mathcal{S})\tau 
   =  \zeta \mathcal{G}_1(E_\mathcal{S})e_{\lambda, \mathcal{S}}^\infty(\phi)+ \epsilon.
\end{equation}
Therefore, take $ \mathcal{G}_2(E_\mathcal{S}):= \zeta \mathcal{G}_1(E_\mathcal{S})$. 
\end{proof}

\begin{proof} of Theorem \ref{EQU} $ $

 Let $\phi\in Ham(M,\omega)$, and let $\epsilon\textgreater0$. It follows from 
Lemma \ref{Vanishingamma} and Lemma \ref{Nonvanishingamma}
that,
\begin{equation}
 \|\phi\|_{Hofer} + 
\|\rho_1^\epsilon\|_{Hofer}\|\psi^\epsilon_1\|_{Hofer} \leq (1 +
\|\rho^\epsilon_1\|_{Hofer}) \left(\mathcal{G}_2(E_\mathcal{S})e_{\lambda, \mathcal{S}}^\infty(\phi) + \epsilon\right), 
\end{equation}

where $\rho^\epsilon$ and $\psi^\epsilon$  are 
two elements $\rho^\epsilon$ and $\psi^\epsilon$ 
of $Ham(M,\omega)$  with $\phi = \rho^\epsilon\circ\psi^\epsilon$. 
 The desired inequality follows from the fact that $e_{\lambda, \mathcal{S}}^\infty(\phi)\leq 
2\|\phi\|_{\lambda, \mathcal{S}}^\infty$. Therefore, the corresponding constant 
$ \mathcal{G}(E_\mathcal{S})$  can be chosen  as: 
$\mathcal{G}(E_\mathcal{S}): = 2\mathcal{G}_2(E_\mathcal{S})$. 
\end{proof}

\begin{remark} Assume that 
the positive constant $\mathcal{G}(E_\mathcal{S})$ is finite, and let\\ 
$diam_{Hofer}(Ham(M,\omega))$ denote 
the diameter of the group $ Ham(M,\omega)$ with respect to the Hofer metric.  
Then, 
from Theorem \ref{EQU}, we see that the constant $\kappa$ in Proposition \ref{cP11} 
can be chosen 
as: $
 \kappa := \left(\mathcal{G}(E_\mathcal{S})\right)^{-1}\left( 1 + diam_{Hofer}(Ham(M,\omega))\right)^{-1}$. 
\end{remark}


 \begin{center}
 {\bf Acknowledgments:} 
 \end{center}
 \begin{center}
I would like to thank the following Professors for helpful comments concerning the first drafts of this note:
\\ Augustin Banyaga, 
Paul Seidel  and Ferdinand Ngakeu.
\end{center}

\end{document}